\documentclass[notitlepage,leqno,10pt]{article}

\usepackage{latexsym}
\usepackage{amsmath}
\usepackage{amsthm}
\usepackage{amsfonts}
\usepackage{amssymb}
\usepackage{graphicx}
\renewcommand{\a }{\alpha }

\renewcommand{\d}{\delta }

\newcommand{\e }{\epsilon }

\renewcommand{\i }{\iota}

\newcommand{\rh }{\rho }

\renewcommand{\th }{\theta }

\newcommand{\intbar}{\mathop{\int\makebox(-13.5,0){\rule[4pt]{.7em}{0.3pt}}%
\kern-6pt}\nolimits}

\newcommand{\be}{\begin{equation}}
\newcommand{\ee}{\end{equation}}
\newenvironment{pf}{\noindent{\sc Proof}.\enspace}{\rule{2mm}{2mm}\medskip}

\newenvironment{pfnb}{\noindent{\bf Proof}}{\rule{2mm}{2mm}\medskip}

\newcommand{\R}{\mathbb{R}}

\newcommand{\Rs}{\mathbb{R}^S}
\newcommand{\Rsd}{{\mathbb{R}^S}^2}
\newcommand{\Rst}{{\mathbb{R}^S}^3}

\newcommand{\Ra}{\mathbb{R}^\alpha}

\newcommand{\N}{\mathbb{N}}

\DeclareMathOperator{\spt}{spt}
\DeclareMathOperator{\diam}{diam}

\begin{document}

\author{Andrea MONDINO$^{1}$}

\date{}

\title{Existence of Integral $m$-Varifolds minimizing $\int |A|^p$ and $\int |H|^p$, $p>m$, in Riemannian Manifolds}

\newtheorem{lem}{Lemma}[section]
\newtheorem{pro}[lem]{Proposition}
\newtheorem{thm}[lem]{Theorem}
\newtheorem{rem}[lem]{Remark}
\newtheorem{cor}[lem]{Corollary}
\newtheorem{df}[lem]{Definition}
\newtheorem{ex}[lem]{Example}
\newtheorem*{Theorem}{Theorem}
\newtheorem*{Lemma}{Lemma}
\newtheorem*{Proposition}{Proposition}
\newtheorem*{claim}{Claim}

\maketitle

\footnotetext[1]{SISSA, Via Bonomea, 265, 34136 Trieste, Italy, E-mail address: mondino@sissa.it, Tel: +39 040 3787-293}
\

\

\begin{center}
\noindent {\sc abstract}. We prove existence of integral rectifiable $m$-dimensional varifolds minimizing functionals of the type $\int |H|^p$ and $\int |A|^p$ in a given Riemannian $n$-dimensional manifold $(N,g)$, $2\leq m<n$ and $p>m$, under suitable assumptions on $N$ (in the end of the paper we give many examples of such ambient manifolds). To this aim we introduce the following new tools: some monotonicity formulas for varifolds in $\Rs$ involving  $\int |H|^p$, to avoid degeneracy of the minimizer, and a sort of isoperimetric inequality to bound the mass in terms of the mentioned functionals.
\bigskip\bigskip

\noindent{\it Key Words:} 
Curvature varifolds, direct methods in the calculus of variation, geometric measure theory, isoperimetric inequality, monotonicity formula

\bigskip

\centerline{\bf AMS subject classification: }
49Q20, 58E99, 53A10, 49Q05.
\end{center}

\section{Introduction}\label{s:in}
Given an ambient Riemannian manifold $(N,g)$ of dimension $n\geq 3$ (with or without boundary), a classical problem in  differential geometry is to find the smooth immersed $m$-dimensional submanifolds, $2\leq m \leq n-1$, with null mean curvature vector, $H=0$, or with null second fundamental form, $A=0$, namely the minimal (respectively, the totally geodesic) submanifolds of $N$ (for more details about the existence see  Example \ref{ex:Ric>}, Example \ref{ex:MxR}, Theorem \ref{thm:ExNRs}, Theorem \ref{thm:ExRs}, Remark \ref{rem:ExStatVar} and Remark \ref{rem:A=0}).

In more generality, it is interesting to study the minimization problems associated to integral functionals depending on the curvatures of the type
\begin{equation}\label{def:Ep}
E_{H,m}^p(M):=\int_M |H|^p \quad \text{or} \quad E_{A,m}^p(M):=\int_M |A|^p, \quad p\geq 1
\end{equation}
where $M$ is a smooth immersed $m$-dimensional submanifold with mean curvature $H$ and second fundamental form $A$; of course the integrals are computed with respect to the $m$-dimensional measure of $N$ induced on $M$.
A global minimizer, if it exists, of $E_{H,m}^p$ (respectively of $E_{A,m}^p$) can be seen as a generalized minimal (respectively totally geodesic) $m$-dimensional submanifold in a natural integral sense.

An important example of such functionals is given by the Willmore functional for surfaces $E_{H,2}^2$ introduced by Willmore (see \cite{Will}) and studied in the euclidean space (see for instance the works of Simon \cite{SiL}, Kuwert and Sch\"atzle \cite{KS}, Rivi\`ere \cite{Riv}) or in Riemannian manifolds (see, for example, \cite{LM},  \cite{Mon1} and \cite{Mon2}). 

The general integral functionals \eqref{def:Ep} depending on the curvatures of immersed submanifolds have been studied, among others, by  Allard \cite{Al}, Anzellotti-Serapioni-Tamanini \cite{AST}, Delladio \cite{Del},  Hutchinson \cite{Hu1}, \cite{Hu2}, \cite{Hu3}, Mantegazza \cite{MantCVB} and Moser \cite{Mos}.

In order to get the existence of a minimizer, the technique adopted in the present paper (as well as in most of the aforementioned ones) is the so called direct method in the calculus of variations. As usual, it is necessary to enlarge the space where the functional is defined and to work out a compactness-lowersemicontinuity theory in the enlarged domain.  

In the present work, the enlarged domain is made of generalized $m$-dimensional submanifolds of the fixed ambient Riemannian manifold $(N,g)$: the integral rectifiable $m$-varifolds introduced by Almgren in \cite{Alm} and by Allard in \cite{Al}. Using integration by parts formulas, Allard \cite{Al} and Hutchinson \cite{Hu1}-Mantegazza \cite{MantCVB} defined a weak notion of mean curvature and of second fundamental form respectively (for more details about this part see Section \ref{AppVar}).  
Moreover these objects have good compactness and lower semicontinuity properties with respect to the integral functionals  above.\newline

The goal of this paper is to prove  existence and partial regularity of an $m$-dimensional  minimizer (in the enlarged class of the rectifiable integral $m$-varifolds with weak mean curvature or with generalized second fundamental form in the sense explained above) of  functionals of the type \eqref{def:Ep}. Actually we will consider more general functionals modeled on this example, see Definition \ref{def:F} for the expression of the considered integrand $F$.

More precisely, given a compact subset $N \subset \subset \bar{N}$ of an $n$-dimensional Riemannian manifold $(\bar{N},g)$ (which, by Nash Embedding Theorem, can be assumed isometrically embedded in some $\Rs$) we will denote 
\begin{eqnarray}
HV_m(N)&:=& \{V\text{ integral rectifiable $m$-varifold of $N$ with weak mean curvature $H^N$ relative to $\bar{N}$} \}\nonumber \\
CV_m(N)&:=& \{V\text{ integral rectifiable $m$-varifold of $N$ with generalized second fundamental form $A$} \};\nonumber
\end{eqnarray}
for more details see Section \ref{AppVar}; in any case, as written above, the non expert reader can think about the elements of $HV_m(N)$ (respectively of $CV_m(N)$) as generalized $m$-dimensional submanifolds with mean curvature $H^N$ (respectively with second fundamental form $A$).
Precisely, we consider the following two minimization problems

\begin{equation}\label{def:beta}
\beta_{N,F}^m:=\inf \left\{ \int_{G_m(N)} F(x,P,H^N) dV: V\in HV_m(N), V\neq 0 \text{ with weak mean curvature } H^N \text{ relative to } \bar{N} \right\}
\end{equation}
and
\begin{equation}\label{def:alpha}
\alpha_{N,F}^m:=\inf \left\{ \int_{G_m(N)} F(x,P,A) dV: V\in CV_m(N), V\neq 0 \text{ with generalized second fundamental form } A \right\}
\end{equation}
where $F$ is as in Definition \ref{def:F} and satisfies \eqref{def:FHp} ( respectively \eqref{def:FAp}). As the reader may see, the expressions  $\int_{G_m(N)} F(x,P,H^N) dV$ (respectively $\int_{G_m(N)} F(x,P,A) dV$) are the natural generalizations of the functionals $E_{H,m}^p$ (respectively $E_{A,m}^p$) in \eqref{def:Ep} with $p>m$ in the context of varifolds.
\\Before stating the two main theorems, let us recall that an integral rectifiable $m$-varifold $V$ on $N$ is associated with a ``generalized $m$-dimensional subset'' $\spt \mu_V$ of $N$ together with an integer valued density function $\th(x)\geq 0$ which carries the ``multiplicity'' of each point  (for the precise definitions, as usual, see Section \ref{AppVar}). 

At this point we can state the two main theorems of this work. Let us start with the mean curvature.

\begin{thm}\label{thm:ExMinH}
Let $N\subset \subset \bar{N}$ be a compact subset with non empty interior, $int(N)\neq \emptyset$, of the $n$-dimensional Riemannian manifold $(\bar{N},g)$ isometrically embedded in some $\Rs$ (by Nash Embedding Theorem), fix $m\leq n-1$ and consider a function $F:G_m(N)\times \Rs \to \R^+$   satisfying \eqref{def:F} and \eqref{def:FHp}, namely
$$F(x,P,H)\geq C|H|^p$$
for some $C>0$ and $p>m$.

Then, at least one of the following two statement is true:

a) the space $(N,g)$ contains a non zero $m$-varifold with null weak mean curvature $H^N$ relative to $\bar{N}$ (in other words, $N$ contains a stationary $m$-varifold; see Remark \ref{rem:StatVar} for the details),

b)  the  minimization problem \eqref{def:beta} corresponding to $F$ has a solution i.e. there exists a non null integral $m$-varifold $V\in HV_m(N)$ with weak mean curvature $H^N$ relative to $\bar{N}$ such that 
$$\int_{G_m(N)} F(x,P,H^N) dV=\beta _{N,F} ^m=\inf \left\{ \int_{G_m(N)} F(x,P,\tilde{H}^N) d\tilde{V}: \tilde{V}\in HV_m(N), \tilde{V}\neq 0\right\}. $$
Moreover, in  case $b)$ is true, we have $\beta _{N,F} ^m>0$ and the minimizer $V$ has the following properties:

b1) the varifold $V$ is indecomposable in  $HV_m(N)$

b2) the diameter of $\spt \mu_V$ as a subset of the Riemannian manifold $(\bar{N},g)$ is strictly positive
$$\diam_{\bar{N}} (\spt \mu_V)>0. $$
\end{thm}

For the precise notion of indecomposability see Definition \ref{def:decomposable}, intuitively (for more details see Remark \ref{rem:decomposable}) the statement $b1)$ is telling that the support $\spt \mu_V$ of the spatial measure $\mu_V$ associated to $V$ is connected.
\begin{rem}\label{rem:HighRegH}
It could be interesting to study the regularity of the minimizer $V$. Since $H \in L^p(V)$, $p>m$ given by \eqref{def:FHp}, we have the following structure result of Allard:  fixed $x\in \spt \mu_V$, under the hypothesis that the density in $x$ satisfies $\th(x)=1$ plus other technical assumptions (see Theorem 8.19 in \cite{Al}),  $\spt \mu_V $ is locally a graph of a $C^{1,1-\frac{m}{p}}$ function. Moreover, under similar assumptions, Duggan proved local $W^{2,p}$ regularity in \cite{Dug}. In the multiple density case the regularity problem is more difficult. For instance,  in \cite{Brak}, is given an example of a varifold $\tilde{V}$ with bounded weak mean curvature whose spatial support contains a set $C$ of strictly positive measure such that if  $x \in C$ then $\spt \mu_{\tilde{V}}$ does not correspond to the graph of even a multiple-valued function in any neighborhood of $x$.  Nevertheless, recently, Menne proved that an integral varifold with locally bounded first variation is $C^2$-rectifiable (see Theorem 1 in \cite{MenJGA}). 
\end{rem}

Now let us state the second main Theorem about the second fundamental form $A$.

\begin{thm}\label{thm:ExMin}
Let $N\subset \subset \bar{N}$ be a compact subset with non empty interior, $int(N)\neq \emptyset$, of the $n$-dimensional Riemannian manifold $(\bar{N},g)$ isometrically embedded in some $\Rs$ (by Nash Embedding Theorem), fix $m\leq n-1$ and consider a function $F:G_m(N)\times \Rst\to \R^+$   satisfying \eqref{def:F} and \eqref{def:FAp}, namely
$$F(x,P,A)\geq C|A|^p$$
for some $C>0$ and $p>m$.

Then, at least  one of the following two statements is true:

a) the space $(N,g)$ contains a non zero $m$-varifold with null generalized second fundamental form,

b)  the  minimization problem \eqref{def:alpha} corresponding to $F$ has a solution i.e. there exists a non null curvature $m$-varifold $V\in CV_m(N)$ with generalized second fundamental form $A$ such that 
$$\int_{G_m(N)} F(x,P,A) dV=\a _{N,F} ^m=\inf \left\{ \int_{G_m(N)} F(x,P,\tilde{A}) d\tilde{V}: \tilde{V}\in CV_m(N), \tilde{V}\neq 0\right\}. $$
Moreover, in case $b)$ is true, we have $\a _{N,F} ^m>0$ and the minimizer $V$ has the following properties:

b1) the varifold  $V$ is indecomposable in $CV_m(N)$ (see Definition \ref{def:decomposable}),

b2) the diameter of $\spt \mu_V$ as a subset of the Riemannian manifold $(\bar{N},g)$ is strictly positive
$$\diam_{\bar{N}} (\spt \mu_V)>0, $$

b3) For every $x\in \spt \mu_V$, $V$ has a unique tangent cone at $x$ and this tangent cone is a finite union of $m$-dimensional subspaces $P_i$ with integer multiplicities $m_i$; moreover, in some neighborhood of $x$ we can express $V$ has a finite union of graphs of $C^{1,1-\frac{m}{p}}$, $m_i$-valued functions defined on the respective affine spaces $x+P_i$ ($p$ given in \eqref{def:FAp}).  
\end{thm}

\begin{rem}\label{rem:HighRegA}
For the precise definitions and results concerning  $b3)$, the interested reader can look at the original paper \cite{Hu2} of Hutchinson. Notice that the boundary of $N$ does not create problems since, by our definitions, the minimizer $V$ is a fortiori an integral $m$-varifold with generalized second fundamental form $A \in L^p(V), p>m,$ in the $n$-dimensional Riemannian manifold $(\bar{N},g)$ which has no boundary. Moreover, by Nash Embedding Theorem, we can assume $\bar{N}\subset \Rs$; therefore $V$ can be seen as an integral $m$-varifold with generalized second fundamental form $A \in L^p(V),$ $p>m,$ in $\Rs$ and the regularity theorem of Hutchinson can be applied. 

It could be interesting to prove higher regularity of the minimizer $V$. About this point, notice that it is not trivially true that $V$ is locally a union of graphs of  $W^{2,p}$ (Sobolev) functions. Indeed in \cite{AGP} there is an example of a curvature $m$-varifold $\tilde{V} \in CV_m(\Rs)$, $S\geq 3$, $2\leq m\leq S-1$, with second fundamental form in $L^p$, $p>m$, which is not a union of graphs of $W^{2,p}$  functions.  

In the spirit of proving higher regularity of the minimizer of such functionals we mention the preprint of Moser \cite{Mos} where the author proves smoothness of the minimizers of $\int |A|^2$ in the particular case of codimension $1$ Lipschitz graphs in $\Rs$. 

Recall also the aforementioned result of Menne \cite{MenJGA} giving that a curvature $m$-varifold is $C^2$-rectifiable.
\end{rem}
  
In both theorems, a delicate point is whether or not $a)$ is satisfied (fact which trivializes the result); we will study this problem in Section \ref{Sec:ExRem}: we will recall two general classes of examples (given by White in \cite{Whi}) of Riemannian manifolds with boundary where $a)$ is not satisfied in codimension $1$, we will give two new examples for higher codimensions (namely Theorem \ref{thm:ExNRs} and Theorem \ref{thm:ExRs}) and we will propose a related open problem in Remark \ref{rem:A=0}. Here, let us just remark that every compact subset $N\subset \subset \Rs$ for $s>1$ does not satisfy $a)$ (see Theorem \ref{thm:ExRs}). \newline

The idea for proving the results is to consider a minimizing sequence $\{V_k\}_{k\in\N}$ of varifolds, show that it is compact (i.e. there exists a varifold $V$ and a subsequence $\{V_{k'}\}$ converging to $V$ in an appropriate sense) and it is non degenerating: if the masses decrease to $0$ the limit would be the null varifold so not a minimizer, and if the diameters decrease to $0$ the limit would be a point which has no geometric relevance.  

In order to perform the analysis of the minimizing sequences, in Section \ref{Sec:MF} we prove monotonicity formulas for integral rectifiable $m$-varifolds in $\Rs$ with weak mean curvature in $L^p$, $p>m$. These formulas are similar in spirit to the ones obtained by Simon in \cite{SiL} for smooth surfaces in $\Rs$  involving the Willmore functional. These estimates are a fundamental tool for proving the non degeneracy of the minimizing sequences and we think they might have other applications.

To show the compactness of the minimizing sequences it is crucial to have a uniform upper bound on the masses   (for the non expert reader: on the volumes of the generalized submanifolds). Inspired by the paper of White \cite{Whi}, in Section \ref{Sec:IsoIne} we prove some isoperimetric inequalities involving our integral functionals which give the mass bound on the minimizing sequences in case $a)$ in the main theorems is not satisfied. The compactness follows and is proved in the same Section. Also in this case, we think that the  results may have other interesting applications.

The proofs of the two main theorems is contained in Section \ref{Sec:min} and \ref{Sec:minH}. Finally, as written above, Section \ref{Sec:ExRem} is devoted to examples and remarks: we will notice that a large class of manifold with boundary can be seen as compact subset of manifold without boundary, we will give examples where the assumption for the isoperimetric inequalities are satisfied and we will end with a related open problem.

The new features of the present paper relies, besides the main theorems, in the new tools introduced in Section \ref{Sec:MF} and Section \ref{Sec:IsoIne}, and in the new examples presented in Section \ref{Sec:ExRem}.  
 
\begin{center}

{\bf Acknowledgments}

\end{center}

\noindent 
This work has been supported by M.U.R.S.T under the Project FIRB-IDEAS ``Analysis and Beyond''.
\\The author would like to thank  G. Bellettini, E. Kuwert, A. Malchiodi, C. Mantegazza, U. Menne and B. White for stimulating and fundamental discussions about the topics of this paper.

\section{Notations, conventions and basic concepts on varifolds} \label{AppVar}
First of all, large positive constants are always denoted by $C$, and the value of
$C$ is allowed to vary from formula to formula and also within the
same line. When we want to stress the dependence of the constants on
some parameter (or parameters), we add subscripts to $C$, as $C_N$,
etc.. Also constants with subscripts are allowed to vary.\newline

In this Section we review the concept of curvature varifold introduced by Hutchinson in \cite{Hu1} giving a slightly more general definition; namely Hutchinson defines the curvature varifolds as ``special''  integral varifolds in a Riemannian manifold but, as a matter of facts, the same definition makes sense for an even non rectifiable varifold in a subset of a Riemannian manifold. So we will define (a priori non rectifiable) varifolds with curvature, which are endowed with a generalized second fundamental form.

We start by recalling some basic facts about varifolds. For more details, the interested reader may look at the standard references \cite{Fed}, \cite{Mor}, \cite{SiGMT} or, for  faster introductions, at \cite{Mant} or the appendix of \cite{Whi}.
 
Consider a (maybe non compact) $n$-dimensional Riemannian manifold $(\bar{N},g)$. Without loss of generality, by the Nash Theorem, we can assume that 
$$(\bar{N},g) \hookrightarrow \Rs  \text{ isometrically embedded for some } S>0.$$ 
We will be concerned with a  subset $N\subset \bar{N}$  which, a fortiori, is also embedded in $\Rs$: $N\hookrightarrow \Rs$. Since throughout the paper $N \subset \subset \bar{N}$ is a compact subset (in the end of the article we will also assume  that it has non empty interior $int(N) \neq \emptyset$) also in this Section is assumed to be so, even if most of the following definitions and properties are valid for more general subsets.

Let us denote with $G(S,m)$ the Grassmaniann of unoriented $m$-dimensional linear subspaces of $\Rs$, with
$$G_m(\bar{N}):=(\Rs \times G(S,m)) \cap \{(x,P):x\in \bar{N}, P\subset T_x\bar{N} \text{ $m$-dimensional linear subspace} \}$$
and with
$$G_m(N):=G_m(\bar{N}) \cap \{(x,P):x\in N\}. $$

 We recall that an $m$-\emph{varifold} $V$ on $N$ is a Radon measure on $G_m(N)$ and that the sequence of varifolds $\{V_k\}_{k\in \N}$ converges to the varifold $V$ in varifold sense if $V_k \to V$ weak as Radon measures on $G_m(N)$; i.e. 
$$\int_{G_m(N)}  \phi \; dV_k\to \int_{G_m(N)}  \phi \; dV$$
as $k\to \infty$, for all $\phi \in C^0_c(G_m(N))$. 
A special class of varifolds are the \emph{rectifiable varifolds}: given  a countably $m$-rectifiable, ${\cal H}^m$- measurable subset $M$ of $N\subset \Rs$ and $\th$ a non negative locally ${\cal H} ^m$ integrable function on $M$, the \emph{rectifiable varifold} $V$ associated to $M$ and $\th$ is defined as
$$V(\phi):=\int_{M} \th(x) \phi(x,T_xM) d {\cal H}^m\quad \forall \phi\in C^0_c(G_m(N))$$
and sometimes is denoted with $V(M,\th)$. Recall that if $M,\th$ are as above then the approximate tangent space $T_xM$ exists for ${\cal H}^m$-almost every $x\in M$ (Theorem 11.6 in \cite{SiGMT}, for the definitions see 11.4 of the same book). If moreover $\th$ is integer valued, then we say that $V$ is an \emph{integral varifold}; the set of the integral $m$-varifolds in $N$ is denoted by $IV_m(N)$.

If $V$ is a $k$-varifold, let $|V|$ denote its mass:
$$|V|:=V(G_m(N)).$$
Observe that we have a natural projection 
\begin{equation}\label{eq:pi}
\pi:G_m(N)\to N \quad  (x,P)\mapsto x,
\end{equation}
and   pushing forward the measure $V$ via the projection $\pi$, we have a positive Radon measure $\mu_V$ on $N$ 
$$\mu_V(B):=V(\pi^{-1}(B))=V(G_m(B)) \quad \forall B\subset N \text{ Borel set}.$$
Since $V$ is a measure on $G_m(N)$, its support is a closed subset of $G_m(N)$. If we project that closed set on $N$ by the projection $\pi$ then we get the \emph{spatial support} of $V$, which coincides with $\spt \mu_V$. 
 
Now let us define the notion of measure-function pair.
 
\begin{df} \label{def:MFpair}
Let $V$ be a Radon measure on $G_m(N)$ (i.e. a varifold) and $f:G_m(N) \to \Ra$ be a well defined $V$ almost everywhere $L^1_{loc}(V)$ function. Then we say that $(V,f)$ is a \emph{measure-function pair} over $G_m(N)$ with values in $\Ra$.

Given $\{(V_k,f_k)\}_{k\in \N}$ and $(V,f)$ measure-function pairs over $G_m(N)$ with values in $\Ra$, suppose $V_k  \to V$ weak as Radon measures in $G_m(N)$ (or equivalently as varifolds in $N$). Then we say \emph{$(V_k,f_k)$ converges to $(V,f)$ in the weak sense} and write 
$$(V_k,f_k)\rightharpoonup (V,f)$$
if   $V_k \lfloor f_k \to V\lfloor f$ weak convergence of Radon vector valued measures. In other words, if
$$\int_{G_m(N)} \left\langle f_k, \phi \right\rangle dV_k\to \int_{G_m(N)} \left\langle f, \phi \right\rangle dV$$
as $k\to \infty$, for all $\phi \in C^0_c(G_m(N), \Ra)$, where $\left\langle .,. \right\rangle$ is the scalar product in $\Ra$. 
\end{df} 

\begin{df} \label{def:F}
Suppose $F:G_m(N)\times \Ra \to \R$. We will denote the variables in $G_m(N)\times \Ra$ by (x,P,q). We say that $F$ satisfies the condition \eqref{def:F} if the following statements are verified:
\\i) $F$ is continuous,
\\ii) $F$ is non negative ( i.e. $F(x,P,q)\geq 0$ for all $(x,P,q) \in G_m(N)\times \Ra$) and  $F(x,P,q)=0$ if and only if $q=0$,
\\iii) $F$ is convex in the $q$ variables, i.e.
$$F(x,P,\lambda q_1 +(1-\lambda) q_2) \leq \lambda F(x,P,q_1)+ (1-\lambda) F(x,P,q_2)$$
for all $\lambda \in (0,1), (x,P) \in G_m(N), q_1 \in \Ra, q_2\in \Ra$,
\\iv) $F$ has superlinear growth in the $q$ variables, i.e. there exists a continuous function $\phi$, where $\phi:G_m(N)\times [0,\infty)\to [0,\infty)$, $0\leq \phi(x,P,s) \leq \phi (x,P,t)$ for $0\leq s \leq t$ and $(x,P)\in G_m(N)$, $\phi(x,P,t)\to \infty$ locally uniformly in $(x,P)$ as $t\to \infty$, such that 
$$\phi(x,P,|q|)|q|\leq F(x,P,q)$$
for all $(x,P,q)\in G_m(N)\times \Ra$.  
\end{df}

An example (trivial but fundamental for this paper) of such an $F$ is $F(x,P,q):=|q|^p$ for any $p>1$.

\begin{rem}\label{rem:LSCF}
For simplicity, in Definition \ref{def:F}, we assumed the same conditions of Hutchinson (\cite{Hu1} Definition 4.1.2) on $F$ but some hypotheses can be relaxed. For example, about the results in this paper, if $F=F(q)$ depends only on the $q$ variables it is enough to assume (in place of i)) that $F$ is lower semicontinuous (see Theorem 6.1 in \cite{MantCVB}).
\end{rem}
In the aforementioned paper, Hutchinson proves the following useful compactness and lower semicontinuity Theorem (see Theorem 4.4.2 in \cite{Hu1}):

\begin{thm} \label{thm:CompLsc}
Suppose $\{(V_k,f_k)\}_{k\in \N}$ are measure-function pairs over $G_m(N)$ with values in $\Ra$. Suppose $V$ is a Radon measure on $G_m(N)$ (i.e a varifold in $N$) and $V_k \to V$ weak converges as Radon measures (equivalently  varifold converges in $N$). Suppose $F:G_m(N)\times \Ra \to \R$ satisfies the condition \eqref{def:F}. Then the following are true:
\\i) If there exists $C>0$ such that 
\begin{equation}\label{eq:Fbounded} 
\int_{G_m(N)} F(x,P,f_k(x,P)) dV_k\leq C \quad \forall k\in \N
\end{equation}
then there exists a function $f \in L^1_{loc}(V)$ such that, up to subsequences, $(V_k,f_k)\rightharpoonup (V,f)$.
\\ii)  if there exists $C>0$ such that \eqref{eq:Fbounded} is satisfied and $(V_k,f_k)\rightharpoonup (V,f)$ then 
$$\int_{G_m(N)} F(x,P,f(x,P)) dV \leq \liminf_k  \int_{G_m(N)} F(x,P,f_k(x,P)) dV_k.$$
\end{thm}

Now we want to define the varifolds of $N$ with curvature. Observe that given $(x,P)\in G_m(N)$, the $m$-dimensional linear subspace $P\subset T_x \bar{N}\subset \Rs$ can be identified with the orthogonal projection matrix on $Hom(\Rs,\Rs)\cong \Rsd$
$$P\equiv [P_{ij}]\in \Rsd.$$
Similarly, the tangent space of $\bar{N}$ at $x$ can be identified with its orthogonal projection matrix 
$$T_x\bar{N}\equiv Q(x):= [Q_{ij}(x)]\in \Rsd.$$
Before defining the varifolds with curvature let us introduce a bit of notation: given $\phi = \phi (x,P) \in C^1(\Rs \times \Rsd)$ we call the partial derivatives of $\phi$ with respect to the variables $x_i$ and $P_{jk}$ (freezing all other variables) by 
$$D_i\phi \quad \text{and} \quad D^*_{jk}\phi \quad \text{for}\quad i,j,k=1,\ldots,S$$
respectively. In the following definition we will consider the quantity
$$P_{ij}\frac{\partial \psi}{\partial x_j}(x) \quad \text{for } \psi\in C^1(\bar{N});$$
we mean that $\psi$ is extended to a $C^1$ function to some neighborhood of $x\in \Rs$ and, since $P$ is the projection on a $m$-subspace of $T_x\bar{N}$, the definition does not depend on the extension. Observe moreover that the quantity depends on $(x,P)$ so it is a function on $G_m(\bar{N})$.  

\begin{df}\label{def:VC}
Let $V$ be an $m$-varifold on $N\subset  \bar{N}\hookrightarrow \Rs$, $m\leq n-1$. We say that $V$ is a \emph{varifold with (generalized) curvature or with (generalized) second fundamental form} if there exist real-valued functions $B_{ijk}$, for $1\leq i,j,k\leq S$, defined $V$ almost everywhere in $G_m(N)$ such that on setting $B=[B_{ijk}]$ the following are true:
\\i) $(V,B)$ is a measure-function pair over $G_m(N)$ with values in $\Rst$.%, moreover $B$ depends just on the space variable for $V$-a.e. in the sense that
%$$B(x,P)=B(x,O_x(P)) \quad \mu_{V}-a.e.\; x \in N,$$
%for any continuous assignment of orthogonal transformations $O:N\to SO(m), x\mapsto O_x$.  
\\ii) For all functions $\phi = \phi (x,P) \in C^1_c (\Rs \times \Rsd)$ one has 
\begin{equation}\label{eq:VC}
0=\int_{G_m(N)} [P_{ij} \, D_j\phi +B_{ijk} \, D_{jk}^*\phi+B_{jij} \,\phi]dV \quad \text{for } i=1,\ldots,S.
\end{equation}
In this case $B$ is called \emph{(generalized) curvature} and we can also define the \emph{(generalized) second fundamental form } of $V$ (with respect to $\bar{N}$) as the $L^1_{loc}(V)$ function with values in $\Rst$ 
\begin{eqnarray}
A&:&G_m(N) \to \Rst, \nonumber \\
A^{k}_{ij} (x,P)&:=&P_{lj} B_{ikl}(x,P)-P_{lj}  P_{iq} \frac{\partial Q_{kl}}{\partial x_q}(x). \label{def:A}
\end{eqnarray}
We will denote the set of \emph{integral} $m$-varifolds of $N$ with generalized curvature as $CV_m(N)$ and we will call them \emph{curvature $m$-varifolds}. 
A varifold $V$ is said to have \emph{null generalized second fundamental form} if  there exists a generalized curvature $B$ whose associated generalized second fundamental form $A$ is null.
\end{df}

Observe that we use different notation of \cite{Hu1}: we call $B$ what Hutchinson calls $A$ and vice versa; this is because we want to denote with $A$ the second fundamental form with respect to $\bar{N}$. 
%Notice also that, even if not explicitly written, also in the classical definition of Hutchinson (Definition 5.2.1 in \cite{Hu1}) it is assumed that the generalized curvature $B$ depends $V$-a.e. just on the space variable $x$; this is because the varifold $V$ was assumed to be integral, so  admitting a unique approximate tangent space at $\mu_V$-a.e. $x \in N$ where the varifold measure $V$ concentrates. 
Moreover, as  shown in Section 5 of \cite{Hu1}, if $V$ is the integral varifold associated to a smooth immersed $m$-submanifold of $N$ then $A$ coincides with the classical second fundamental form with respect to $N$.

\begin{rem} \label{rem:BA}
By definition, the generalized second fundamental form $A$ is expressed in terms of $B$ but, as Hutchinson proved in \cite{Hu1} Propositions 5.2.4 and 5.2.6, it is possible to express $B$ in terms of $A$. Indeed, choosing appropriate test functions, with some easy computations one can prove that
\begin{equation}\label{eq:BA}
B_{ijk}=A^k_{ij}+A^j_{ik}+ P_{jl} P_{iq} \frac{\partial Q_{lk}}{\partial x_q}(x)+P_{kl} P_{iq} \frac{\partial Q_{lj}}{\partial x_q}(x).
\end{equation}
\end{rem}

\begin{rem}\label{rem:NonUniq}
As U. Menne pointed out to the author in a personal communication, the proof of Hutchinson of the uniqueness of the generalized curvature $B$ (see \cite{Hu1}, Proposition 5.2.2) uses in a crucial way  that, in his case, the varifold is  rectifiable.

At the level of generality of Definition \ref{def:VC}, i.e. $V$ is an a priori non rectifiable varifold, the arguments of Hutchinson does not ensure the uniqueness of the generalized second fundamental form, and to our knowledge the uniquness is open in the general case. 

However throughout the paper we mostly work with integral varifolds, for which the uniqueness of the second fundamental form is granted by \cite{Hu1}; the only place where we need a priori non rectifiable varifolds with generalized second fundamental form is in Theorem  \ref{thm:AreaEst} and its applications, where we use the concept of maybe non-rectifiable varifold with null generalized second fundamental form. We will mean a varifold $V$ such that there exists a generalized curvature $B$ whose associated generalized second fundamental form $A$, as in \eqref{def:A}, is null.
\end{rem}

Now let us recall the fundamental compactness and lower semi continuity Theorem of Hutchinson (Theorem 5.3.2 in \cite{Hu1})
\begin{thm}\label{thm:HuLSC}
Consider $\{V_k\}_{k\in \N} \subset CV_m(N)$ with generalized second fundamental forms $\{A_k\}_{k\in \N}$, $V$ an integral $m$-varifold of $N$ and suppose $V_k \to V$ in varifold sense. Let $F:G_m(N)\times \Rst \to \R$ be a function satisfying the condition \eqref{def:F} and assume that 
$$\int_{G_m(N)} F(x,P,A_k) dV_k \leq C$$
for some $C>0$ independent on $k$. Then

i) $V \in CV_m(N)$ with generalized second fundamental form $A$, 

ii) $(V_k, A_k) \rightharpoonup (V,A)$ in the weak sense of measure-function pairs,

iii) $\int_{G_m(N)} F(x,P,A) dV \leq \liminf_k \int_{G_m(N)} F(x,P,A_k) dV_k$. 
\end{thm} 

Now we briefly recall the definition of first variation of an $m$-varifold $V$ in $\Rs$; the original definitions are much more general, here we recall only the facts we need for this paper.

\begin{df}\label{def:1variat}
Let $V$ be an $m$-varifold in $\Rs$ and let $X$ be a $C^1_c(\Rs)$ vector field. We define \emph{first variation} $\d V$ the linear functional on $C^1_c(\Rs)$ vector fields
$$\d V(X):= \int_{G_m(\Rs)} div_P X(x) dV(x,P);$$
where for every $P\in G(S,m)$, 
$$div_P X:= \sum_{i=1}^S \nabla_i ^P X^i= \sum_{i,j=1}^S P_{ij} D_j X^i,$$
where $\nabla^P f=P (\nabla f)$ is the projection on $P$ of the gradient in $\Rs$ of $f$ and $\nabla_i^P:=e_i \cdot \nabla^P$ (where $\{e_i\}_{i=1,\ldots,S}$ is an orthonormal basis of $\Rs$).

$V$ is said to be of \emph{locally bounded first variation} in $\Rs$ if for every relatively compact open $W\subset \subset \Rs$
there exists a constant $C_W < \infty$ such that
$$|\d V(X)|\leq C_W \sup_W |X|$$
for all $X\in C^1_c(\Rs)$ vector fields with support in $W$.
\end{df}

An interesting subclass of varifolds with locally bounded first variation are the varifolds with weak mean curvature.

\begin{df}\label{def:VarWeakH}
 Let $V$ be an $m$-varifold  in $\Rs$ and  $H:G_m(\Rs) \to \Rs$ an $L^1_{loc}(V)$ function (in the previous notation we would say that $(V,H)$ is a measure-function pair on $G_m(\Rs)$ with values in $\Rs$); then we say that $V$ has \emph{weak mean curvature $H$} if for any vector field $X \in C^1_c(\Rs)$ one has
\begin{equation}\label{def:weakHVar}
\d V(X):= \int_{G_m(\Rs)} div_P X(x) dV(x,P)= -\int_{G_m(\Rs)} H \cdot X dV(x,P) \quad .
\end{equation} 
\end{df}

Observe that if $V=V(M,\th)$ is a rectifiable varifold with weak mean curvature $H$ then with abuse of notation we can write $H(x)=H(x,T_xM)$ and we get the following identities:
\begin{equation}\label{def:weakH}
\int_M div_M X d\mu_V=\int_{G_m(\Rs)} div_{T_xM} X(x) dV=-\int_{G_m(\Rs)} H(x,T_xM) \cdot X dV=-\int_M H(x) \cdot X d\mu_V,
\end{equation}
where $div_M X$ is the tangential divergence of the vector field $X$ and is defined to be $div_M X(x):= div_{T_xM} X(x)$ where $T_xM$ is the approximate tangent space to $M$ at $x$ (which exists for $\mu_V$-a.e. $x$).

\begin{rem}\label{rem:AH}
As Hutchinson observed in \cite{Hu1}, if $V$ is an $m$-varifold on $N\hookrightarrow \Rs$ with generalized curvature $B=[B_{ijk}]_{i,j,k=1,\ldots,S}$ then, as a varifold in $\Rs$, $V$ has weak mean curvature $H_i=\sum_{j=1}^S B_{jij}$ for $i=1,\ldots, S$. Indeed, for any relatively compact open subset  $W\subset \subset \Rs$ and any vector field $X\in C^1_c(\Rs)$ with compact support in $W$, taking $\phi=X^i,\; i=1,\ldots, S$ in equation  \eqref{eq:VC} and summing over $i$ we get 
$$0=\int_{G_m(\Rs)} [P_{ij}D_j X^i(x)+B_{jij}(x,P) X^i(x)]dV(x,P)$$
which implies
$$ \d V(X):= \int_{G_m(\Rs)} div_P X(x) dV(x,P)=- \int_{G_m(\Rs)} B_{jij}(x,P) X^i(x)\; dV(x,P);$$
the conclusion follows from the fact that $B \in L^1_{loc} (V)$.
\end{rem}

Now let us define the varifolds with weak mean curvature in a compact subset $N\subset \subset \bar{N}$ of a Riemannian manifold $(N,g)$ isometrically embedded in $\Rs$. 

\begin{df}\label{def:VH}
Let $V$ be an $m$-varifold on $N\subset  \bar{N}\hookrightarrow \Rs$, $m\leq n-1$. We say that $V$ is a \emph{varifold with weak mean curvature $H^N$ relative to $\bar{N}$} if it has weak mean curvature $H^{\Rs}$ as varifold in $\Rs$. In this case the value of $(H^N)_i,$ $i=1,\ldots,S$ is given by
\begin{equation}\label{eq:weakH}
(H^N)_i=(H^{\Rs})_i- P_{jk} \frac{\partial Q_{ij}}{\partial x^k}.
\end{equation}
Consistently with the notation introduced for the curvature varifolds, we denote with $HV_m(N)$ the set of \emph{integral} $m$-varifolds on $N$ with weak mean curvature relative to $\bar{N}$; the elements of $HV_m(N)$ are called \emph{mean curvature varifolds}.
\end{df}

Observe that in case $V$ is the varifold associated to a smooth submanifold of $\bar{N}$ then $H^N$ coincides with the classical mean curvature relative to $\bar{N}$ (it is enough to trace the identity (i) of Proposition 5.1.1 in \cite{Hu1} recalling that we denote with $A$,$Q$ what Hutchinson calls $B$,$S$). Moreover, as an exercise, the reader may check that also in the general case the vector $\left(P_{jk} \frac{\partial Q_{ij}}{\partial x^k}\right)_{i=1,\ldots,S}$ of $\Rs$  is orthogonal to $\bar{N}$ (fix a point $x$ of $\bar{N}$ and choose a base of $T_x\bar{N}$ in which the Christoffel symbols of $\bar{N}$ vanish at $x$; write down the orthogonal projection matrix $Q$ with respect to this base and check the orthogonality condition).

\begin{rem}\label{rem:dVH}
If $V$ is an $m$-varifold on $N\subset  \bar{N}\hookrightarrow \Rs$, $m\leq n-1$ with weak mean curvature $H^N$ relative to $\bar{N}$ then, for each compactly supported vector field $X \in C^1_c(\bar{N})$ tangent to $\bar{N}$,
$$\d V(X)=\int_{G_m(N)} div_P X(x) dV(x,P)=-\int_{G_m(N)} H^N\cdot X dV(x,P)\quad .$$
\end{rem}

This fact gives consistency to Definition \ref{def:VH} and follows from Definition \ref{def:VH}, from formula \eqref{def:weakHVar} and the orthogonality of  $\left(P_{jk} \frac{\partial Q_{ij}}{\partial x^k}\right)_{i=1,\ldots,S}$ to $\bar{N}$.

\begin{rem}\label{rem:StatVar}
If $V$ is an $m$-varifold on $N\subset \bar{N} \hookrightarrow \Rs$, $m\leq n-1$ with null weak mean curvature $H^N=0$ relative to $\bar{N}$ then, for each compactly supported vector field $X \in C^1_c(\bar{N})$ tangent to $\bar{N}$,
$$\d V(X)=\int_{G_m(N)} div_P X(x) dV(x,P)=0.$$
In this case we say that $V$ is an $m$-varifold in $N$ with null weak mean curvature relative to $\bar{N}$ or, using more classical language, that $V$ is a \emph{stationary} $m$-varifold in $N$ (where stationary as to be intended in $\bar{N}$). 
\end{rem}

\begin{df}\label{def:decomposable}
A curvature  $m$-varifold $V \in CV_m(N)$ (resp. $V\in HV_m(N)$) is said to be \emph{decomposable in $CV_m(N)$} (resp. in $HV_m(N)$) if there exist two non null curvature varifolds $0\neq V_1,V_2 \in CV_m(N)$ (resp. $0\neq V_1,V_2 \in HV_m(N)$) such that $V=V_1+V_2$. If $V$ is not decomposable in $CV_m(N)$ (resp. $HV_m(N)$) it is said \emph{indecomposable in $CV_m(N)$} (resp. in $HV_m(N)$).
\end{df}

\begin{rem}\label{rem:decomposable}
Notice that, given $V\in CV_m(N)$ (resp.  $V \in HV_m(N)$), if the support of the spatial measure $\spt \mu_V$ has two connected components at positive distance then $V$ is decomposable; indeed, localizing Definition \ref{def:VC} using cutoff functions, it is clear that the two connected components detect two non null elements of $CV_m(N)$ whose sum is the original varifold $V$. 
\\In case instead $\spt \mu_V$ has countably many connected components which are accumulating, it is not so clear if each connected component detect an element of $CV_m(N)$ (resp. $HV_m(N)$) as it is not possible to isolate each component by using  cutoff functions. 
\\On the other hand even if $\spt \mu_V$ is connected the varifold $V$ may be decomposable, as the example of two smooth embedded compact  submanifolds with non empty intersection shows. 
\end{rem}

\section{Monotonicity formulas for integral $m$-varifolds with weak mean curvature in $L^p$, $p>m$ }\label{Sec:MF}
Let $V=V(M,\th)$ be an integral varifold of $\Rs$   (associated to the rectifiable set $M\subset \Rs$ and with integer multiplicity function $\th$) with weak mean curvature $H$ (since throughout this section we consider only varifolds in $\Rs$ and there is no ambiguity, we adopt the easier notation $H$ for $H^{\Rs}$). Let us write $\mu$ for $\mu_V:= \pi_{\sharp}(V)$ the push forward of the varifold measure $V$ on $G_m(N)$ to $N$ via the standard projection $\pi:G_m(N)\to N, \pi(x,P)=x$ (see Section \ref{AppVar} for more details); of course $\mu_V$ can also be seen as $\mu_V={\cal H}^m \lfloor \th$, the restriction of the $m$-dimensional Hausdorff measure to the multiplicity function $\th$.

The first Lemma is a known fact (see for example the book of Leon Simon \cite{SiGMT} at page 82) of which we report also the proof for completeness.

\begin{lem}\label{lem:DiffMF}
Let $V=V(M,\th)\in IV_m(\Rs)$ be with weak mean curvature $H$ as above and fix a point $x_0\in M$. For $\mu$-a.e. $x \in M$ call $r(x):=|x-x_0|$ and $D^\perp r$ the orthogonal projection of the gradient vector $Dr$ onto $(T_xM)^\perp$. Consider a nonnegative function $\phi \in C^1(\R)$ such that
$$\phi'(t)\leq0 \; \forall t\in \R, \quad \phi(t)=1 \; \text{for } t\leq \frac{1}{2}, \quad \phi(t)=0 \;\text{for } t\geq 1.$$
For all $\rh>0$ let us denote  
\begin{eqnarray}
I(\rh)&:=&\int_M \phi(r/\rh) d\mu, \nonumber \\
L(\rh)&:=& \int_M \phi(r/\rh) (x-x_0)\cdot H d\mu, \nonumber \\
J(\rh)&:=& \int_M \phi(r/\rh) |D^\perp r|^2 d\mu; \nonumber
\end{eqnarray}
then
\begin{equation}\label{eq:DiffMF}
\frac{d}{d\rh}[\rh^{-m} I(\rh)]=\rh^{-m} J'(\rh)+\rh^{-m-1}L(\rh).
\end{equation}
\end{lem} 

\begin{pf}
The idea is to use formula \eqref{def:weakH} and choose the vector field $X$ in an appropriate way in order to get informations about $V$. First of all let us recall that for any function $f\in C^1(\Rs)$ and any $x \in M$ where the approximate tangent space $T_xM$ exists (it exists for $\mu$-a.e. $x\in M$ see \cite{SiGMT} 11.4-11.6 ) one can define the tangential gradient as the projection of the gradient in $\Rs$ onto $T_xM$:
$$\nabla^M f:= \sum_{j,l=1}^{S} P^{jl} D_l f(x) e_j$$
where $D_l f$ denotes the partial derivative $\frac{\partial f}{\partial x^l}$ of $f$, $P^{jl}$ is the matrix of the orthogonal projection of $\Rs$ onto $T_xM$ and $\{e_j\}_{j=1,\ldots,S}$ is an orthonormal basis of $\Rs$.
Denoted $\nabla_j^M:=e_j \cdot \nabla^M$, recall that the tangential divergence is defined as
$$div_M X:=\sum_{j=1}^S \nabla_j^M X^j;$$
moreover it is easy to check the Leibniz formula
$$div_M fX:= \nabla^M f \cdot X+ f\, div_M X \quad \forall f\in C^1(\Rs) \text{ and } \forall X \in C^1(\Rs)\text{ vector field. }$$ 
Now let us choose the vector field. Fix $\rh>0$ and consider the function $\gamma \in C^1(\R)$ defined as
$$\gamma(t):=\phi(t/\rh);$$
then of course we have the following properties: 
$$\gamma'(t)\leq0 \; \forall t\in \R, \quad \gamma(t)=1 \; \text{for } t\leq \frac{\rh}{2}, \quad \gamma(t)=0 \;\text{for } t\geq \rh.$$
Call $r(x):=|x-x_0|$ and choose the vector field
$$X(x):= \gamma(r(x)) (x-x_0).$$
Using the Leibniz formula we get
\begin{eqnarray}
div_M X&=&\nabla^M \gamma(r) \cdot (x-x_0)+ \gamma(r) div_M (x-x_0) \nonumber\\
&=& r \gamma'(r) \frac{(x-x_0)^T}{|x-x_0|} \frac{(x-x_0)^T} {|x-x_0|}+ m \gamma(r) \nonumber\\
&=& r \gamma'(r) (1-|D^\perp r|^2)+ m \gamma(r),\label{eq:divMX}
\end{eqnarray}
where $u^T$ is the projection of the vector $u \in \Rs$ onto $T_pM$ and  $D^\perp r= \frac{(x-x_0)^\perp}{|x-x_0|}$ is the orthogonal projection of the gradient vector $Dr$ onto $(T_xM)^\perp$.
The equation \eqref{def:weakH} of the weak mean curvature thus yields
\begin{equation}\label{eq:gamma}
m \int_M \gamma(r) d\mu + \int_M r \gamma'(r) d\mu=\int_M r \gamma'(r) |D^\perp r|^2 d\mu -\int_M H\cdot (x-x_0) \gamma(r) d\mu.
\end{equation}
Now recall that $\gamma(r)=\phi(r/\rho)$, so $r\gamma'(r)=\frac{r}{\rh} \phi'(r/\rh)=-\rh \frac{\partial }{ \partial \rh} [\phi(r/\rh)].$
Thus, combining \eqref{eq:gamma} and  the definitions of $I(\rh),J(\rh)$ and $L(\rh)$ one gets 
$$m I(\rh)-\rh I'(\rh)=-\rh J'(\rh)-L(\rh).$$
Thus, multiplying both sides by $\rh^{-m-1}$ and rearranging we obtain
$$\frac{d}{d\rh}[\rh^{-m}I(\rh)]=\rh^{-m} J'(\rh)+\rh^{-m-1} L(\rh).$$ 
This concludes the proof
\end{pf}

Estimating from below the right hand side of \eqref{eq:DiffMF} and integrating, we get the following useful inequalities.

\begin{pro}\label{pro:LocMF}
Let $V=V(M,\th)\in IV_m(\Rs)$ be with weak mean curvature $H\in L^p(V)$, $p>m$ (we mean that $\int_{G_m(\Rs)} |H|^p dV<\infty$ or equivalently, denoted with an abuse of notation $H(x)=H(x,T_xM)$, $\int_M |H|^p d\mu<\infty$). Fixed a point $x_0\in M$ and $0<\sigma<\rh<\infty$, then 
\begin{equation}
[\sigma^{-m} \mu(B_\sigma(x_0))]^{\frac{1}{p}}\leq [\rh^{-m} \mu(B_\rh(x_0))]^{\frac{1}{p}}+\frac{p^2}{p-m} \rh^{1-\frac{m}{p}} \Big( \int_{B_\rh(x_0)} |H|^p d\mu \Big)^{\frac{1}{p}}- \frac{p^2}{p-m} \sigma^{1-\frac{m}{p}} \Big( \int_{B_\sigma(x_0)} |H|^p d\mu \Big)^{\frac{1}{p}}.
\end{equation}
\end{pro}

\begin{pf}
Let us estimate from below the right hand side of equation \eqref{eq:DiffMF}. Observe that 
$$J'(\rh)=\frac{d}{d\rh} \int_M \phi(r/\rh) |D^{\perp}r|^2 d\mu=-\rh^{-2} \int_M r \phi'(r/\rh)|D^{\perp}r|^2 d\mu\geq 0 $$
since $\phi'(t)\leq 0$ for all $t\in \R$. 
Thus we can say that
\begin{equation}\label{eq:MFprov}
\frac{d}{d\rh}[\rh^{-m} I(\rh)]\geq \rh^{-m-1}L(\rh).
\end{equation}
Let us estimate from below the right hand side by Schwartz inequality:
\begin{eqnarray}
\rh^{-m-1}L(\rh)&=&\rh^{-m-1} \int_M \phi(r/\rh) (x-x_0)\cdot H d\mu \nonumber \\
                &\geq& - \rh^{-m-1} \int_M \big(\phi(r/\rh)^{\frac{1}{p}} |H|\big)\; |x-x_0| \phi(r/\rh)^{\frac{p-1}{p}}  d\mu. \nonumber
\end{eqnarray}
Now recalling that $\phi(t)=0$ for $t\geq 1$ we get that $\phi(r/\rh)=0$ for $r\geq \rh$ so $|x-x_0|$ in the integral can be estimated from above by $\rh$ and we can say that
$$ \rh^{-m-1}L(\rh)\geq - \rh^{-m} \int_M \big(\phi(r/\rh)^{\frac{1}{p}} |H|\big)\;  \phi(r/\rh)^{\frac{p-1}{p}}  d\mu; $$
thus, by Holder inequality, for all $p>1$
\begin{eqnarray}
\rh^{-m-1}L(\rh)&\geq&- \rh^{-m} \Big(\int_M \phi(r/\rh) |H|^p d\mu\Big)^{\frac{1}{p}} \Big(\int_M \phi(r/\rh) d\mu \Big)^{\frac{p-1}{p}} \nonumber\\
&=& - \rh^{-m} \Big(\int_M \phi(r/\rh) |H|^p d\mu\Big)^{\frac{1}{p}}\; I(\rh)^{\frac{p-1}{p}}\label{eq:L}.
\end{eqnarray}
Putting together inequalities \eqref{eq:MFprov} and \eqref{eq:L} we get
$$\frac{d}{d\rh}[\rh^{-m} I(\rh)]\geq - \rh^{-m} \Big(\int_M \phi(r/\rh) |H|^p d\mu\Big)^{\frac{1}{p}}\; I(\rh)^{\frac{p-1}{p}};$$
multiplying both sides by $\rh^{m-\frac{m}{p}} I(\rh)^{\frac{1}{p}-1}$ and rearranging we get
$$\frac{d}{d\rh}[\rh^{-m} I(\rh)]^{\frac{1}{p}}\geq-p \;\rh^{-\frac{m}{p}} \Big(\int_M \phi(r/\rh) |H|^p d\mu\Big)^{\frac{1}{p}}.$$
Now, after choosing $p>m$, integrate the last inequality from $\sigma$ to $\rh$ (the same $\rh$ chosen in the statement of the Proposition) and get with an integration by parts of the right hand side
\begin{align}
\rh^{-\frac{m}{p}}I(\rh)^{\frac{1}{p}}-\sigma^{-\frac{m}{p}} I(\sigma)^{\frac{1}{p}} &\geq -p\int_{\sigma} ^{\rh} \Big[ \left(t^{-\frac{m}{p}}\right) \Big( \int_M\phi(r/t) |H|^p d\mu \Big)^{\frac{1}{p}}  \Big]dt \nonumber\\
&=-p\Big[\Big(1-\frac{m}{p}\Big)^{-1} \Big(\rh^{1-\frac{m}{p}} \Big(\int_{M}\phi(r/\rh) |H|^p d\mu \Big)^{\frac{1}{p}}-\sigma^{1-\frac{m}{p}} \Big(\int_{M}\phi(r/\sigma) |H|^p d\mu \Big)^{\frac{1}{p}} \Big)\Big]\nonumber\\
&\qquad +p \int_\sigma ^\rh \Big[ \Big(1-\frac{m}{p}\Big)^{-1} t^{1-\frac{m}{p}} \Big(\frac{d}{dt} \int_M \phi(r/t) |H|^p d\mu \Big)\Big] dt \label{eq:provI}
\end{align}
Observe that, as before  for $J'(\rh)$, since $\phi'(t)\leq 0$ for all $t$ it follows
$$\frac{d}{dt} \int_M \phi(r/t) |H|^p d\mu=-t^{-2} \int_M r\phi'(r/t) |H|^p d\mu \geq 0$$ 
so the second integral in equation \eqref{eq:provI} is non negative and, recalling the definition of $I$, we can write
\begin{eqnarray}\label{eq:provMF}
\Big(\rh^{-m}\int_M \phi(r/\rh) d\mu \Big)^{\frac{1}{p}}- \Big(\sigma^{-m}\int_M \phi(r/\sigma) d\mu \Big)^{\frac{1}{p}}&\geq& \frac{p^2}{p-m}  \Big[-\rh^{1-\frac{m}{p}} \Big(\int_{M}\phi(r/\rh) |H|^p d\mu \Big)^{\frac{1}{p}} \nonumber \\
&& \qquad \qquad+ \sigma^{1-\frac{m}{p}} \Big(\int_{M}\phi(r/\sigma) |H|^p d\mu \Big)^{\frac{1}{p}}\Big].
\end{eqnarray}
Now observe that during all this proof and during all the proof of Lemma \ref{lem:DiffMF} the only used properties of $\phi$ have been
$$\phi\in C^1(\R), \quad \phi'(t)\leq 0 \;\forall t\in \R,\quad \phi(t)\leq 1 \; \forall t\in \R,\quad \phi(t)=0 \;\forall t\geq 1; $$
thus, for all such $\phi$, the inequality \eqref{eq:provMF} holds.
Now taking a sequence $\phi_k$ of such functions pointwise converging to the characteristic function of $]-\infty,1]$ and, using the Dominated Convergence Theorem, passing to the limit on $k$ in \eqref{eq:provMF}  we get
$$\big[\rh^{-m} \mu (B_\rh(x_0)) \big]^{\frac{1}{p}}- \big[\sigma^{-m} \mu (B_\sigma(x_0))\big]^{\frac{1}{p}}\geq \frac{p^2}{p-m}  \Big[-\rh^{1-\frac{m}{p}} \Big(\int_{B_\rh(x_0)}|H|^p d\mu \Big)^{\frac{1}{p}}+ \sigma^{1-\frac{m}{p}} \Big(\int_{B_\sigma(x_0)}|H|^p d\mu \Big)^{\frac{1}{p}}\Big]. $$
Rearranging we can conclude that
$$\big[\sigma^{-m} \mu (B_\sigma(x_0))\big]^{\frac{1}{p}}\leq \big[\rh^{-m} \mu (B_\rh(x_0)) \big]^{\frac{1}{p}}+ \frac{p^2}{p-m}  \rh^{1-\frac{m}{p}} \Big(\int_{B_\rh(x_0)}|H|^p d\mu \Big)^{\frac{1}{p}}- \frac{p^2}{p-m}\sigma^{1-\frac{m}{p}} \Big(\int_{B_\sigma(x_0)}|H|^p d\mu \Big)^{\frac{1}{p}}$$
\nopagebreak[4]
\end{pf}

From Corollary 17.8 page 86 of \cite{SiGMT}, if $H\in L^p(V)$ for some $p>m$, then the density $\th(x)=\lim_{\rh\downarrow 0}\frac{\mu(\bar{B}_\rh(x))}{w_m \rh^m}$ exists at every point $x\in \Rs$ and $\th$ is an upper semicontinuous function. Hence, letting $\sigma \to 0$, one has
$$[\omega_m \th(x_0)]^{\frac{1}{p}}\leq \Big[\frac{\mu(B_\rh(x_0))}{\rh^m} \Big]^{\frac{1}{p}}+\frac{p^2}{p-m} \Big[\rh^{p-m} \int_{B_\rh(x_0)} |H|^p d\mu \Big]^{\frac{1}{p}}. $$
Using the inequality $a^{\frac{1}{p}}+b^{\frac{1}{p}}\leq 2^{\frac{p-1}{p}} (a+b)^{\frac{1}{p}}$ given by the concavity of the function $t\mapsto t^{\frac{1}{p}}$ with $p>1$ and $t>0$, we get 
$$\omega_m \th(x_0)\leq 2^{p-1} \Big[\frac{\mu(B_\rh(x_0))}{\rh^m} +\Big(\frac{p^2}{p-m}\Big)^p \rh^{p-m} \int_{B_\rh(x_0)} |H|^p d\mu \Big].$$
Since $V\in IV_m(\Rs)$, then $\th$ is integer valued and by definition $\th\geq 1$ $\mu$-a.e. From the upper semicontinuity of $\th$ it follows that $\th(x)\geq 1$ for all $x\in \spt \mu$ (where, as before, $\mu$ is the spatial measure associated to $V$). Then the last  formula can be written more simply getting the $\mathbf{fundamental\; inequality}$ 
  
\begin{equation}\label{eq:FundIn}
1\leq C_{p,m} \Big[\frac{\mu(B_\rh(x_0))}{\rh^m} + \rh^{p-m} \int_{B_\rh(x_0)} |H|^p d\mu \Big] \quad \forall x_0\in \spt \mu,
\end{equation}
where $C_{p,m}>0$ is a positive constant depending on $p,m$ and such that $C_{p,m}\to \infty$ if $p\downarrow m$.

\begin{rem}
The fundamental inequality can be extended to the case $p=m$ by using the isoperimetric inequality, see 8.3 in \cite{Al};  for discussions on related results see Proposition 3.1 in \cite{LeMa} and  2.5 in \cite{MenACV}.
\end{rem}

Using the fundamental inequality now we can link through inequalities  the mass of $V$,  the diameter of $M$ and the $L^p$ norm of the weak mean curvature $H$.

\begin{lem}\label{A<dH}
Let $V=V(M,\th)\in IV_m(\Rs)$ be a non null integral $m$-varifold with compact spatial support $\spt \mu \subset \Rs$ and weak mean curvature $H\in L^p(V)$ for some $p>m$. Then, called $d=\diam_{\Rs}(\spt \mu)$ the diameter of $\spt \mu$ as a subset of $\Rs$,
\begin{equation}\label{eq:A<dH}
|V|\leq \left(\frac{d}{m}\right)^p  \int_M |H|^p d\mu.
\end{equation}
\end{lem}

\begin{pf}
In the same spirit of the proof of Lemma \ref{lem:DiffMF} we choose a suitable vector field $X$ to plug in the mean curvature equation \eqref{def:weakH} $$\int_M div_M X d\mu =-\int_M X\cdot H d\mu $$
in order to get informations about the varifold $V=V(M,\th)$.
 Now fix a point $x_0 \in \spt \mu$ and  simply let $X(x)=x-x_0$. Since  $div_M X=m$ $\mu$-a.e. (for more details see the proof of Lemma \ref{lem:DiffMF}), observing that $|X|\leq d$ $\mu$-a.e. and estimating the right hand side by Holder inequality we get 
$$m |V|\leq d \Big(\int_M |H|^p d\mu \Big)^{\frac{1}{p}} |V|^{\frac{p-1}{p}}. $$
Now multiply both sides by $|V|^{\frac{1}{p}-1}$ and raise to the power $p$ in order to get the thesis.  
 
\end{pf}

\begin{lem}\label{d<AH}
Let $V=V(M,\th)\in IV_m(\Rs)$ be a non null integral indecomposable (in $HV_m(\Rs)$) $m$-varifold with compact spatial support $\spt \mu\subset \Rs$ and weak mean curvature $H\in L^p(V)$ for some $p>m$. Then, called $d=\diam_{\Rs}(\spt \mu)$,
\begin{equation}\label{eq:d<AH}
d \leq C_{p,m} \Big(\int_M |H|^p d\mu \Big)^{\frac{m-1}{p}} |V|^{1-\frac{m-1}{p}}
\end{equation}
where $C_{p,m}>0$ is a positive constant depending on $p,m$ and such that $C_{p,m}\to \infty$ if $p\downarrow m$.
\end{lem} 

\begin{pf}
Since $\spt \mu \subset \Rs$ is compact, then there exist $x_0,y_0\in \spt \mu$ such that 
$$d=|x_0-y_0|.$$
Let $\rh \in ]0,d/2]$ and call $N:=\lfloor d/\rh\rfloor$ the integer part of $d/\rh$. For every $j=1,\ldots, N-1$ the indecomposability of $V$ implies the existence of $r_j\in \left(j+\frac{1}{4}, j+\frac{3}{4}\right)$ such that
$$B_{r_j \rh}(y_0)\cap \spt \mu \neq \emptyset. $$ 
Let $y_j\in \partial B_{r_j \rh}(y_0)\cap \spt \mu$ and observe that for each ball $B_{\rh/4}(y_j), \; j=0,\ldots, N-1$ we have the fundamental inequality \eqref{eq:FundIn}; since the balls  $B_{\rh/4}(y_j), \; j=0,\ldots, N-1$ are pairwise disjoint, summing up over $j$ we get
\begin{equation*}
N\leq C_{p,m} \left( \frac{|V|}{\rh^m} + \rh^{p-m} \int_M |H|^p d\mu \right).
\end{equation*}
Moreover, since $N=\lfloor d/\rh\rfloor\geq \frac{d}{2\rh}$,  we have
\begin{equation}\label{eq:FundInN}
d\leq 2\rh N \leq C_{p,m} \left( \frac{|V|}{\rh^{m-1}} + \rh^{p-m+1} \int_M |H|^p d\mu \right). 
\end{equation}
Now let us choose $\rh$ in an appropriate way; observe that taken
$$\rh=\frac{m}{2} \left(\frac{|V|}{\int_M |H|^p d\mu}\right)^{\frac{1}{p}},$$
in force of the estimate \eqref{eq:A<dH}, the condition $\rh \leq d/2$ is satisfied. Finally, plugging this value of $\rh$ into equation \eqref{eq:FundInN}, after some trivial computation we conclude that
$$d\leq C_{p,m} \;|V|^{\frac{p-m+1}{p}} \left(\int_M |H|^p d\mu \right)^{\frac{m-1}{p}}. $$
\end{pf}

\begin{rem}
Notice that, in Lemma \ref{d<AH}, the assumption that $V$ is indecomposable in $HV_m(\Rs)$ can be replaced by asking that the support of the spatial measure $\spt \mu$ is connected.
\end{rem}

Combining the Fundamental Inequality with the previous lemmas we are in position to prove a lower diameter and mass bound. 
\begin{lem}\label{lem:d>H}
Let $V=V(M,\th)\in IV_m(\Rs)$ be a non null  integral $m$-varifold with  spatial support $\spt \mu \subset \Rs$ and weak mean curvature $H\in L^p(V)$ for some $p>m$. Then, called $d:=\diam_{\Rs}(\spt \mu)$
\begin{equation}\label{eq:d>H}
d\geq \frac{1}{C_{p,m}\Big(\int_M |H|^p d\mu \Big)^{\frac{1}{p-m}}}
\end{equation}
where $C_{p,m}>0$ is a positive constant depending on $p,m$ and such that $C_{p,m}\to \infty$ if $p\downarrow m$.
\end{lem}

\begin{pf}
If $d=\infty$, the inequality \eqref{eq:d>H} is trivially satisfied; hence we can assume that $\spt \mu \subset \Rs$ is compact. It follows that there exist $x_0,y_0\in \spt \mu$ such that 
$$d=|x_0-y_0|.$$
Recall the Fundamental Inequality \eqref{eq:FundIn} and choose $\rh=d$ obtaining
\begin{equation}\label{eq:1VH}
1\leq C_{p,m} \Big(\frac{|V|}{d^m} + d^{p-m} \int_{M} |H|^p d\mu \Big).
\end{equation}
From Lemma \ref{A<dH}, 
$$|V|\leq \frac{1}{m^p} d^p  \int_M |H^p| d\mu, $$
hence the inequality \eqref{eq:1VH} becomes
$$1\leq C_{p,m}\;  d^{p-m} \int_{M} |H|^p d\mu $$
and we can conclude.

\end{pf}

\begin{lem}\label{A>H}
Let $V=V(M,\th)\in IV_m(\Rs)$ be a non null integral $m$-varifold with compact spatial support $\spt \mu \subset \Rs$ and weak mean curvature $H\in L^p(V)$ for some $p>m$. Then
\begin{equation}\label{eq:A>H}
|V| \geq  \frac{1}{C_{p,m}\Big(\int_M |H|^p d\mu \Big)^{\frac{m}{p-m}}}
\end{equation}
where $C_{p,m}>0$ is a positive constant depending on $p,m$ and such that $C_{p,m}\to \infty$ if $p\downarrow m$.
\end{lem} 

\begin{pf}
First of all if $V$ is decomposable in $HV_m(\Rs)$ then, by definition, each component is  an integral varifold with weak mean curvature in $L^p$. Hence we can assume that $V$ is indecomposable in $HV_m(\Rs)$, otherwise just argue on a non null component and observe that the inequality \eqref{eq:A>H} is well behaved.

Call as before $d:=\diam_{\Rs} (\spt \mu)$; from  the inequality \eqref{eq:d<AH},
\begin{equation*}
|V|\geq \frac{d^{\frac{p}{p-m+1}}}{ \left(\int_M |H|^p d\mu\right)^{\frac{m-1}{p-m+1}}}.
\end{equation*}
But from the last inequality \eqref{eq:d>H}, 
\begin{equation*}
d^{\frac{p}{p-m+1}} \geq \frac{1}{C_{p,m}\Big(\int_M |H|^p d\mu \Big)^{\frac{p}{(p-m)(p-m+1)}}}.
\end{equation*}
Combining the two estimates, with an easy computation we get the conclusion.

\end{pf}

\begin{pro} \label {pro:HausConv}
Let $\{V_k=V_k(M_k,\th_k)\}_{k\in \N} \subset IV_m(\Rs)$ be a sequence of integral varifolds with weak mean curvature $H_k\in L^p(V_k)$  for some $p>m$ and associated spatial measures $\mu_k$. Assume a uniform bound on the $L^p$ norms of $H_k$:
$$\exists C>0:\; \forall k\in \N\quad \int_{M_k} |H_k|^p d\mu_k= \int_{G_m(\Rs)} |H_k|^p dV_k\leq C, $$
and assume a uniform bound on the spatial supports $\spt \mu_k$:
$$\exists R>0: \spt \mu_k\subset B_R^{\Rs} $$
where $B_R^{\Rs}$ is the ball of radius $R$ centered in the origin in $\Rs$.
 
It follows that if there exists a Radon measure $\mu$ on $\Rs$ such that 
$$\mu_k \to \mu\quad \text{weak as Radon measures,}$$
then
$$\spt \mu_k \to \spt \mu \quad \text{in Hausdorff distance sense.}$$
 
\end{pro}

\begin{pf}
First of all observe that the uniform bound on the spatial supports $\spt \mu_k$ implies that $\spt \mu $ is compact. Since $\spt \mu$ is compact, recall that $\spt \mu_k \to \spt\mu$ if and only if the set of the all possible limit points of all possible sequences $\{x_k\}_{k\in \N}$ with $x_k \in \spt \mu_k$ coincides with $\spt \mu$. Let us prove it by double inclusion.

i) since $\mu_k \to \mu$ weak as Radon measures of course $\forall x \in spt\,\mu$ there exists a sequence $\{x_k\}_{k\in \N}$ with $x_k \in \spt \mu_k$ such that $x_k\to x$. Otherwise there would exist $\e>0$ such that for infinitely many $k'$
$$B_\e(x)\cap \spt \mu_{k'} = \emptyset.$$
This would imply that $\mu_{k'}(B_\e(x))=0$, but $x\in \spt \mu$ so we reach the contradiction 
$$0<\mu(B_\e(x))=\lim_{k'} \mu_{k'}B_\e(x)=0.$$

ii) Let $\{x_k\}_{k\in \N}$ with $x_k \in \spt \mu_k$ be such that $x_k\to x$. We have to show that $x\in \spt \mu$. Let us argue by contradiction:
\\if $x \notin \spt \mu$ then there exists $\e_0>0$ such that
\begin{equation}\label{eq:xsptmu}
0=\mu(B_{\e_0}(x))=\lim_k \mu_k(B_{\e_0}(x)).
\end{equation}
Since $\spt \mu_k \ni x_k\to x$, then for every $\e \in (0,\e_0/2)$ there exists $K_\e>0$ large enough such that
$$x_k\in (\spt \mu_k\cap B_\e(x)) \quad \forall k>K_\e.$$
Now consider the balls $B_\e(x_k)$ for $k>K_\e$: by the triangle inequality $B_\e(x_k)\subset B_{\e_0}(x)$, moreover, since by construction $x_k\in \spt \mu_k$, we can apply the fundamental inequality \eqref{eq:FundIn} to each $B_\e(x_k)$ and obtain
\begin{eqnarray}
1&\leq& C_{p,m} \Big[\frac{\mu_k(B_\e(x_k))}{\e^m} + \e^{p-m} \int_{B_\e(x_k)} |H_k|^p d\mu_k \Big]\nonumber \\
 &\leq& C_{p,m} \Big[\frac{\mu_k(B_{\e_0}(x))}{\e^m} + \e^{p-m} \int_{M_k} |H_k|^p d\mu_k \Big] \quad \forall k>K_\e\label{eq:FundInmuk}.
\end{eqnarray}
Keeping in mind \eqref{eq:xsptmu}, for every fixed $\e \in (0,\e_0/2)$ we can pass to the limit on $k$ in inequality \eqref{eq:FundInmuk} and get
$$\liminf_k \int_{M_k} |H_k|^p d\mu_k\geq \frac{1}{C_{p,m}\, \e^{p-m}}. $$
But $\e>0$ can be arbitrarily small, contradicting the uniform bound $\int_{M_k} |H_k|^p d\mu_k\leq C$ of the assumptions.  
 
\end{pf}

\section{Isoperimetric inequalities and  compactness results}\label{Sec:IsoIne}
\subsection{An isoperimetric inequality involving the generalized second fundamental form} \label{SubSec:IsoIneA}

The following Isoperimetric Inequality involving the generalized second fundamental form is inspired by the paper of White \cite{Whi} and uses the concept of varifold with second fundamental form introduced by Hutchinson \cite{Hu1}. Actually we need a slight generalization of the definition of curvature varifold given by Hutchinson: in Definition 5.2.1 of \cite{Hu1}, the author considers only integral varifolds but, as a matter of facts, a similar definition makes sense for a general varifold. In  Section \ref{AppVar} we  recalled  the needed concepts.

\begin{thm}\label{thm:AreaEst}
Let $N \subset \subset \bar{N}$ be  a compact subset of  a (maybe non compact) $n$-dimensional Riemannian manifold $(\bar{N},g)$ (which, by Nash Embedding Theorem we can assume isometrically embedded in some $\Rs$) and let $m\leq n-1$. Then the following  conditions are equivalent:

i) $N$ contains no nonzero $m$-varifold with null generalized second fundamental form

ii) There is an increasing function $\Phi:\R^+ \to \R^+$ with $\Phi(0)=0$ and a function $F:G_m(N)\times \Rst \to \R^+ $ satisfying \eqref{def:F} (see Section \ref{AppVar}) such that for every $m$-varifold $V$ in $N$ with generalized second fundamental form $A$
$$ |V| \leq \Phi \Big(\int_{G_m(N)}F(x,P,A(x,P)) dV\Big).$$

iii) for every function $F:G_m(N)\times \Rst \to \R^+ $ satisfying \eqref{def:F} (see Section \ref{AppVar}) there exists a constant $C_F>0$ such that for every $m$-varifold $V$ in $N$ with generalized second fundamental form $A$
$$ |V| \leq C_F \int_{G_m(N)}F(x,P,A(x,P)) dV.$$

\end{thm}

\begin{pf}
Of course iii) $\Rightarrow$ ii) $\Rightarrow$ i). It remains to prove that i) $\Rightarrow$ iii).
Let us argue by contradiction: assume that iii) is not satisfied and prove that also i) cannot be satisfied.
\\First fix the function $F$. If iii) is not satisfied then there exists a sequence  $\{(V_k,A_k)\}_{k\in \N}$ of $m$-varifolds in $N$ with generalized second fundamental form (see Definition \ref {def:VC}) such that
$$|V_k|\geq k \int_{G_m(N)} F(x,P,A_k(x,P)) dV_k.$$
We can assume that $|V_k|=1$ otherwise replace $V_k$ with the normalized varifold $\tilde{V}_k:=\frac{1}{|V_k|}V_k$ (observe that the second fundamental form is invariant under this rescaling of the measure and that $\int_{G_m(N)} F(x,P,A_k) dV_k= |V_k|\int_{G_m(N)} F(x,P,A_k) d\tilde{V}_k$). Hence
$$\int_{G_m(N)} F(x,P,A_k(x,P)) dV_k\leq \frac{1}{k}.$$
Recall that $|V_k|=1$ so, from Banach-Alaoglu and Riesz Theorems, there exists a varifold $V$ such that, up to subsequences, $V_k \to V$ in varifold sense (i.e weak convergence of Radon measures on $G_m(N)$). Of course $|V|=\lim_k |V_k|=1$.

Using the notation of \cite{Hu1} (see the Section \ref{AppVar}) we have that the measure-function pairs $(V_k,A_k)$ over  $G_m(N)$, up to subsequences, satisfy the assumptions of Theorem \ref{thm:CompLsc}. From (i) of the mentioned Theorem \ref{thm:CompLsc}, it follows that there exists a measure-function pair $(V,\tilde{A})$ with values in $\Rst$ (i.e a Radon measure $V$ on $G_m(N)$ and a matrix valued function $\tilde{A}\in L^1_{loc}(V)$ ) such that $(V_k,A_k) \rightharpoonup (V,\tilde{A})$ (i.e $V_k \lfloor A_k \to V\lfloor \tilde{A}$ weak convergence of Radon vector valued measures). 

From Remark \ref{rem:BA}  we can express the generalized curvatures $B_k$ of the varifolds $V_k$ in terms of the second fundamental forms $A_k$. Moreover, calling $B$ the corresponding quantity to $\tilde{A}$, from the explicit expression \eqref{eq:BA} it is clear that the weak convergence $(V_k,A_k)\rightharpoonup (V,\tilde{A})$ implies the weak convergence $(V_k,B_k)\rightharpoonup (V,B)$.
\\Passing to the limit in $k$ in \eqref{eq:VC} we see that $(V,B)$ satisfies the equation, so $V$ is an $m$-varifold with generalized curvature $B$. %(the fact that $B$ depends $V$-a.e. just on the space variable is a consequence of Lemma \ref{lem:fk(x)} together with the observation that $B_k$ depend $V_k$-a.e. just on the space variable; the last fact being granted by the assumption that the varifolds $V_k$ are rectifiable and hence admit $\mu_{V_k}$-a.e. a unique approximate tangent plane on which the varifold measure concentrates). 

Now let us check that the corresponding generalized second fundamental form (in sense of equation \eqref{def:A}) to $B$  is $\tilde{A}$. 
\\Call $$\Lambda^{l}_{ij} (x,P):=P_{pj} B_{ilp}(x,P)-P_{pj}  P_{iq} \frac{\partial Q_{lp}}{\partial x_q}(x)$$
the corresponding second fundamental form to $B$ and $\Lambda_k=A_k$ the corresponding to $B_k$.

Since $(V_k,B_k)\rightharpoonup (V,B)$, from the definitions it is clear that $(V_k,\Lambda_k)\rightharpoonup (V,\Lambda)$; but, from the definition of $\tilde{A}$, $(V_k,\Lambda_k)=(V_k,A_k)\rightharpoonup (V,\tilde{A})$. It follows that $\Lambda=\tilde{A}\;$ $V$-almost everywhere and that $\tilde{A}$ is the generalized second fundamental form of $V$ associated to $B$.  

Finally, the lower semicontinuity of the functional (sentence (ii) of Theorem \ref{thm:CompLsc}) implies 
$$\int_{G_m(N)} F(x,P,\tilde{A}) dV\leq \liminf_{k} \int_{G_m(N)} F(x,P,A_k) dV_k=0.$$
From the assumption ii) of condition \eqref{def:F} on $F$ it follows that $\tilde{A}=0$ $V$-almost everywhere; henceforth we constructed a non null $m$-varifold $V$ in $N$ with null second fundamental form and this concludes the proof.
\end{pf}

\begin{rem}
A trivial but fundamental example of $F:G_m(N)\times \Rst \to \R$ satisfying the assumptions of Theorem \ref{thm:AreaEst} is $F(x,P,A)=|A|^p$ for any $p>1$. Hence the Theorem implies that if a compact subset $N$ of a  Riemannian $n$-dimensional manifold $(\bar{N},g)$ does not contain any non null $k$-varifold  ($k\leq n-1$) with null generalized second fundamental form then for every $p>1$ there exists a constant $C_p>0$ such that
$$|V|\leq C_p \int_{G_m(N)}|A|^p dV $$
for every $k$-varifold $V$ in $N$ with generalized second fundamental form $A$.
\end{rem}

Putting together the fundamental compactness and lower semicontinuity Theorem \ref{thm:HuLSC} of Hutchinson  and the Isoperimetric Theorem \ref{thm:AreaEst} we get the following useful compactness-lower semicontinuity result.

\begin{thm}\label{thm:CompLSC} 
Let $N \subset \subset \bar{N}$ be  a compact subset of  a (maybe non compact) $n$-dimensional Riemannian manifold $(\bar{N},g)$ (which, by Nash Embedding Theorem we can assume isometrically embedded in some $\Rs$), fix $m\leq n-1$  and let $F:G_m(N)\times \Rst\to \R^+$  be a function satisfying \eqref{def:F}. 

Assume that, for some $m\leq n-1$, the space $(N,g)$ does not contain any non zero $m$-varifold with null generalized second fundamental form.
  
Consider a sequence $\{V_k\}_{k\in \N} \subset CV_m(N)$ of curvature varifolds  with generalized second fundamental forms $\{A_k\}_{k\in \N}$ such that 
$$\int_{G_m(N)} F(x,P,A_k) dV_k \leq C$$
for some $C>0$ independent on $k$. 

Then there exists  $V \in CV_m(N)$ with generalized second fundamental form $A$ such that, up to subsequences, 

i) $(V_k, A_k) \rightharpoonup (V,A)$ in the weak sense of measure-function pairs,

ii) $\int_{G_m(N)} F(x,P,A) dV \leq \liminf_k \int_{G_m(N)} F(x,P,A_k) dV_k$. 
\end{thm}

\begin{pf}
From Theorem \ref{thm:AreaEst} there exists a constant $C_F>0$ depending on the function $F$ such that $|V_k| \leq C_F \int_{G_m(N)}F(x,P,A_k(x,P)) dV_k$, thus from the boundness of $\int_{G_m(N)} F(x,P,A_k) dV_k $ we have the uniform mass bound
\begin{equation}\label{eq:MB}
|V_k|\leq C
\end{equation}
for some $C>0$ independent on $k$. This mass bound, together with Banach Alaoglu and Riesz Theorems, implies that there exists an $m$-varifold $V$ on $N$ such that, up to subsequences, $V_k \to V$ in varifold sense.

In order to apply Hutchinson compactness Theorem \ref{thm:HuLSC} we have to prove that $V$ actually is an integral $m$-varifold. 
\\From assumption $iv)$ on $F$ of Definition \ref{def:F}, there exists a continuous function  $\phi:G_m(N)\times [0,\infty)\to [0,\infty)$, with $0\leq \phi(x,P,s) \leq \phi (x,P,t)$ for $0\leq s \leq t$ and $(x,P)\in G_m(N)$, $\phi(x,P,t)\to \infty$ locally uniformly in $(x,P)$ as $t\to \infty$, such that 
\begin{equation}\label{eq:phiF}
\phi(x,P,|A|)|A|\leq F(x,P,A)
\end{equation}
for all $(x,P,A)\in G_m(N)\times \Rst.$ 
Since $N$ is compact, also  $G_m(N)$ is so and from the properties of $\phi$ there exists $C>0$ such that $\phi(x,P,|A|)\geq 1$ for $|A|> C$ and any $(x,P)\in G_m(N)$. Thus for every $k$ we can split the computation of the $L^1(V_k)$ norm of $A_k$ as

$$\int_{G_m(N)}|A_k|dV_k=\int_{G_m(N)\cap\{|A_k|\leq C\}}|A_k|dV_k+\int_{G_m(N)\cap\{|A_k|> C\}}|A_k|dV_k. $$
The first term is bounded above by the mass bound \eqref{eq:MB}. About the second term observe that, for $|A|>C$ the inequality \eqref{eq:phiF} implies that $|A|\leq F(x,P,A)$; then also the second term is bounded in virtue on the assumption that $\int_{G_m(N)} F(x,P,A_k) dV_k$ is uniformly bounded.

We have proved that there exists a constant $C$ such that, for all $k\in \N$,
\begin{equation}\label{eq:boundAk}
\int_{G_m(N)}|A_k|dV_k \leq C.
\end{equation}

Now, change point of view and look at the varifolds $V_k$ as curvature varifolds in $\Rs$.
Recall (see Remark \ref{rem:BA}) that the curvature function $B$  can be written in terms of the generalized second fundamental form $A$ relative to $\bar{N}$ and of the extrinsic curvature of the manifold $\bar{N}$ (as submanifold of $\Rs$) which is uniformly bounded on $N$ from the compactness assumption. Using the triangle inequality together with estimate \eqref{eq:boundAk} and the mass bound \eqref{eq:MB} we obtain the uniform estimate of the $L^1(V_k)$ norms of the curvature functions $B_k$
 \begin{equation}\label{eq:boundBkC}
\int_{G_m(\Rs)}|B_k|dV_k \leq C
\end{equation}
for some $C>0$ independent on $k$.

Estimate \eqref{eq:boundBkC} and Remark \ref{rem:AH} tell us that the integral varifolds $V_k$ of $\Rs$ have uniformly bounded first variation: there exists a $C>0$ independent on $k$ such that
$$|\d V_k(X)|\leq C \sup_{\Rs} |X|, \quad \forall X\in C^1_c(\Rs) \text{ vector field}. $$
The uniform bound on the first variations and on the masses of the integral varifolds $V_k$ allow us to apply Allard's integral compactness Theorem (see for example \cite{SiGMT} Remark 42.8 or the original paper of Allard \cite{Al}) and say that the limit varifold $V$ is actually integral. 

The conclusions of the Theorem then follow from Hutchinson Theorem \ref{thm:HuLSC}. 
\end{pf}

\begin{cor}\label{cor:CompMinSeq}
Let $N \subset \subset \bar{N}$ be  a compact subset with non empty interior, $int(N)\neq \emptyset$, of  a (maybe non compact) $n$-dimensional Riemannian manifold $(\bar{N},g)$ (which, by Nash Embedding Theorem can be assumed isometrically embedded in some $\Rs$) and let $F:G_m(N)\times \Rst\to \R^+$  be a function satisfying \eqref{def:F}.  

Assume that, for some $m\leq n-1$, the space $(N,g)$ does not contain any non zero $m$-varifold with null generalized second fundamental form.

Call
\begin{equation}\label{eq:DefAlpha}
\alpha_{N,F}^m:=\inf \left\{ \int_{G_m(N)} F(x,P,A) dV: V\in CV_m(N), V\neq 0 \text{ with generalized second fundamental form } A \right\}
\end{equation} 
and consider a minimizing  sequence $\{V_k\}_{k\in \N} \subset CV_m(N)$ of curvature varifolds  with generalized second fundamental forms $\{A_k\}_{k\in \N}$ such that 
$$\int_{G_m(N)} F(x,P,A_k) dV_k  \downarrow \a _{N,F} ^m.$$

Then there exists  $V \in CV_m(N)$ with generalized second fundamental form $A$ such that, up to subsequences,  

i) $(V_k, A_k) \rightharpoonup (V,A)$ in the weak sense of measure-function pairs,

ii) $\int_{G_m(N)} F(x,P,A) dV \leq \a _F ^m$. 
\end{cor}

\begin{pf}
We only have to check that $\a _{N,F} ^m< \infty$, then the conclusion follows from Theorem \ref{thm:CompLSC}. But the fact is trivial since $int(N)\neq \emptyset$, indeed  we can always construct a smooth compact $m$-dimensional embedded submanifold of $N$, which of course is a curvature  $m$-varifold with finite energy.
\end{pf}

\begin{rem}\label{rem:NonMinA}
Notice that, a priori, Corollary \ref{cor:CompMinSeq} does not ensure the existence of a minimizer since it can happen that the limit $m$-varifold $V$ is null. In the next Section \ref{Sec:min} we will see that, if $F(x,P,A)\geq C |A|^p$ for some $C>0$ and $p>m$, then this is not the case and we have a non trivial minimizer. 
\end{rem}
 
\subsection{An isoperimetric inequality involving the weak mean curvature}\label{SubSec:IsoIneH}

In this Subsection we adapt to the context of varifolds with weak mean curvature the results of the previous Subsection \ref{SubSec:IsoIneA} about varifolds with generalized second fundamental form (for the basic definitions and properties see Section \ref{AppVar}). The following Isoperimetric Inequality involving the weak mean curvature can be seen as a variant of Theorem 2.3 in \cite{Whi}. 

\begin{thm}\label{thm:AreaEstH}
Let $N \subset \subset \bar{N}$ be  a compact subset of  a (maybe non compact) $n$-dimensional Riemannian manifold $(\bar{N},g)$ (which, by Nash Embedding Theorem we can assume isometrically embedded in some $\Rs$) and let $m\leq n-1$. Then the following  conditions are equivalent:

i) $N$ contains no nonzero $m$-varifold with null weak mean curvature relative to $\bar{N}$ (i.e $N$ contains no nonzero stationary $m$-varifold; see Remark \ref{rem:StatVar}).

ii) There is an increasing function $\Phi:\R^+ \to \R^+$ with $\Phi(0)=0$ and a function $F:G_m(N)\times \Rs \to \R^+ $ satisfying \eqref{def:F} (see Section \ref{AppVar}) such that for every $m$-varifold $V$ in $N$ with weak mean curvature $H^N$ relative to $\bar{N}$ 
$$ |V| \leq \Phi \Big(\int_{G_m(N)}F(x,P,H^N(x,P)) dV\Big).$$

iii) for every function $F:G_m(N)\times \Rs \to \R^+ $ satisfying \eqref{def:F} (see Section \ref{AppVar}) there exists a constant $C_F>0$ such that for every $m$-varifold $V$ in $N$  with weak mean curvature $H^N$ relative to $\bar{N}$ 
$$ |V| \leq C_F \int_{G_m(N)}F(x,P,H^N(x,P)) dV.$$

\end{thm}

\begin{pf}
The proof is similar to the proof of Theorem \ref{thm:AreaEst}. Of course iii) $\Rightarrow$ ii) $\Rightarrow$ i). We prove by contradiction that i) $\Rightarrow$ iii): assume that iii) is not satisfied and show that also i) cannot be satisfied.
\\First fix the function $F$. If iii) is not satisfied then there exists a sequence  $\{V_k\}_{k\in \N}$ of $m$-varifolds in $N$ with weak mean curvatures $H^N_k$ relative to $\bar{N}$ (see Definition \ref {def:VH}) such that
$$|V_k|\geq k \int_{G_m(N)} F(x,P,H^N_k(x,P)) dV_k.$$
We can assume that $|V_k|=1$ otherwise replace $V_k$ with the normalized varifold $\tilde{V}_k:=\frac{1}{|V_k|}V_k$ (observe that the weak mean curvature is invariant under this rescaling of the measure and that $\int_{G_m(N)} F(x,P,H^N_k) dV_k= |V_k|\int_{G_m(N)} F(x,P,H^N_k) d\tilde{V}_k$). Hence
$$\int_{G_m(N)} F(x,P,H^N_k(x,P)) dV_k\leq \frac{1}{k}.$$
Recall that $|V_k|=1$ so, from Banach-Alaoglu and Riesz Theorems, there exists a varifold $V$ such that, up to subsequences, $V_k \to V$ in varifold sense (i.e weak convergence of Radon measures on $G_m(N)$). Of course $|V|=\lim_k |V_k|=1$.

Now the measure-function pairs $(V_k,H^N_k)$ over  $G_m(N)$, up to subsequences, satisfy the assumptions of Theorem \ref{thm:CompLsc} and (i) (of the mentioned Theorem \ref{thm:CompLsc}) implies that there exists a measure-function pair $(V,\tilde{H}^N)$ with values in $\Rs$  such that $(V_k,H_k^N) \rightharpoonup (V,\tilde{H}^N)$ weak convergence of measure-function pairs (see Definition \ref{def:MFpair}). 

At this point we have to check that $V$ is an $m$-varifold of $N$ with weak mean curvature $\tilde{H}^N$ relative to $\bar{N}$.
Recall that $N\hookrightarrow \Rs$, so the varifolds $V_k$ can be seen as varifolds with weak mean curvatures $H^{\Rs}_k$ in $\Rs$; from equation \eqref{eq:weakH}, the measure-function pair convergence $(V_k,H^N_k) \rightharpoonup (V,\tilde{H}^N)$ implies the measure-function pair convergence $(V_k,H^{\Rs}_k) \rightharpoonup (V,\tilde{H}^N+P_{jk} \frac{\partial Q_{ij}}{\partial x^k})$ which says ( pass to the limit in Definition \ref{def:VarWeakH}) that $V$ is an $m$-varifold in $\Rs$ with weak mean curvature  $\tilde{H}^N+P_{jk} \frac{\partial Q_{ij}}{\partial x^k}$. Thus, by Definition \ref{def:VH},  $V$ is an $m$-varifold of $N$ with weak mean curvature $H^N:= \tilde{H}^N$ relative to $\bar{N}$.

Finally, the lower semicontinuity of the functional (sentence (ii) of Theorem \ref{thm:CompLsc}) implies 
$$\int_{G_m(N)} F(x,P,H^N) dV\leq \liminf_{k} \int_{G_m(N)} F(x,P,H^N_k) dV_k=0.$$
From the assumption ii) of condition \eqref{def:F} on $F$ it follows that $H^N=0$ $V$-almost everywhere; henceforth we constructed a non null $m$-varifold $V$ in $N$ with null weak mean curvature relative to $\bar{N}$ and this concludes the proof.
\end{pf}

We also have a counterpart of Theorem \ref{thm:CompLSC} concerning the weak mean curvature:
\begin{thm}\label{thm:CompLSCH} 
Let $N \subset \subset \bar{N}$ be  a compact subset of  a (maybe non compact) $n$-dimensional Riemannian manifold $(\bar{N},g)$ (which, by Nash Embedding Theorem we can assume isometrically embedded in some $\Rs$), fix $m\leq n-1$  and let $F:G_m(N)\times \Rs\to \R^+$  be a function satisfying \eqref{def:F}. 

Assume that, for some $m\leq n-1$, the space $(N,g)$ does not contain any non zero $m$-varifold with null weak mean curvature relative to $\bar{N}$.
  
Consider a sequence $\{V_k\}_{k\in \N} \subset HV_m(N)$ of integral $m$-varifolds  with weak mean curvatures $\{H^N_k\}_{k\in \N}$ relative to $\bar{N}$ such that 
$$\int_{G_m(N)} F(x,P,H^N_k) dV_k \leq C$$
for some $C>0$ independent on $k$. 

Then there exists  $V \in HV_m(N)$ integral varifold with weak mean curvature $H^N$ relative to $\bar{N}$ such that, up to subsequences, 

i) $(V_k, H^N_k) \rightharpoonup (V,H^N)$ in the weak sense of measure-function pairs,

ii) $\int_{G_m(N)} F(x,P,H^N) dV \leq \liminf_k \int_{G_m(N)} F(x,P,H^N_k) dV_k$. 
\end{thm}

\begin{pf}
The proof is analogous to the proof of Theorem \ref{thm:CompLSC}. From Theorem \ref{thm:AreaEstH} there exists a constant $C_F>0$ depending on the function $F$ such that $|V_k| \leq C_F \int_{G_m(N)}F(x,P,H^N_k(x,P)) dV_k$, thus from the boudness of $\int_{G_m(N)} F(x,P,H^N_k) dV_k $ we have the uniform mass bound
\begin{equation}\label{eq:MBH}
|V_k|\leq C
\end{equation}
for some $C>0$ independent on $k$. This mass bound, together with Banach Alaoglu and Riesz Theorems, implies that there exists an $m$-varifold $V$ on $N$ such that, up to subsequences, $V_k \to V$ in varifold sense.

The proof that $V$ actually is an integral $m$-varifold is completely analogous to the same statement in the proof of Theorem \ref{thm:CompLSC}: formally substituting $H^N_k$ to $A_k$ in the mentioned proof we arrive to 
\begin{equation}\label{eq:boundHNk}
\int_{G_m(N)}|H^N_k|dV_k \leq C.
\end{equation}

Now, change point of view and look at the varifolds $V_k$ as integral varifolds in $\Rs$.
From Definition \ref{def:VH} the weak mean curvature $H^{\Rs} _k$ in $\Rs$  can be written in terms of $H^N_k$ and of the extrinsic curvature of the manifold $\bar{N}$ (as submanifold of $\Rs$) which is uniformly bounded on $N$ from the compactness assumption. Using the triangle inequality together with estimate \eqref{eq:boundHNk} and the mass bound \eqref{eq:MBH} we obtain the uniform estimate of the $L^1(V_k)$ norms of the weak mean curvatures $H^{\Rs}_k$
 \begin{equation}\label{eq:boundBk}
\int_{G_m(\Rs)}|H^{\Rs}_k|dV_k \leq C
\end{equation}
for some $C>0$ independent on $k$. It follows (see Definition \ref{def:VarWeakH}) that the integral varifolds $V_k$ of $\Rs$ have uniformly bounded first variation: there exists a constant $C>0$ independent on $k$ such that
$$|\d V_k(X)|\leq C \sup_{\Rs} |X|, \quad \forall X\in C^1_c(\Rs) \text{ vector field}. $$
The uniform bound on the first variations and on the masses of the integral varifolds $V_k$ allow us to apply Allard's integral compactness Theorem (see for example \cite{SiGMT} Remark 42.8 or the original paper of Allard \cite{Al}) and say that the limit varifold $V$ is actually integral. 

With the same arguments in the end of the proof of Theorem \ref{thm:AreaEstH}, one can show that the varifold convergence of a subsequence $V_k \to V$ and the uniform energy bound $\int_{G_m(N)} F(x,P,H^N_k) dV_k <C$ implies the existence of a measure-function pair converging subsequence $(V_k,H^N_k) \rightharpoonup (V,H^N)$ for some $\Rs$ -valued function $H^N\in L^1_{loc}(V)$ which actually is the weak mean curvature of $V$ relative to $\bar{N}$. 

We conclude that $V\in HV_m(N)$ is an integral $m$-varifold of $N$ with weak mean curvature $H^N$ relative to $\bar{N}$ and i) holds; property ii) follows from the general Theorem \ref{thm:HuLSC} about measure-function pair convergence (specifically see sentence ii) of the mentioned Theorem). 
\end{pf}

Finally we have a counterpart of Corollary \ref{cor:CompMinSeq}
\begin{cor}\label{cor:CompMinSeqH}
Let $N \subset \subset \bar{N}$ be  a compact subset with non empty interior, $int(N)\neq \emptyset$, of  a (maybe non compact) $n$-dimensional Riemannian manifold $(\bar{N},g)$ (which, by Nash Embedding Theorem can be assumed isometrically embedded in some $\Rs$) and let $F:G_m(N)\times \Rs\to \R^+$  be a function satisfying \eqref{def:F}.  

Assume that, for some $m\leq n-1$, the space $(N,g)$ does not contain any non zero $m$-varifold with null weak mean curvature relative to $\bar{N}$.

Call
\begin{equation}\label{eq:DefBeta}
\beta _{N,F}^m:=\inf \left\{ \int_{G_m(N)} F(x,P,H^N) dV: V\in HV_m(N), V\neq 0 \text{ with weak wean curvature } H^N \text{ relative to } \bar{N} \right\}
\end{equation} 
and consider a minimizing  sequence $\{V_k\}_{k\in \N} \subset HV_m(N)$ of integral varifolds with weak mean curvatures $\{H^N_k\}_{k\in \N}$ such that 
$$\int_{G_m(N)} F(x,P,H^N_k) dV_k  \downarrow \beta _{N,F} ^m.$$

Then there exists an integral $m$-varifold $V \in HV_m(N)$ with weak mean curvature $H^N$ relative to $\bar{N}$ such that, up to subsequences,  

i) $(V_k, H^N_k) \rightharpoonup (V,H^N)$ in the weak sense of measure-function pairs,

ii) $\int_{G_m(N)} F(x,P,H^N) dV \leq \beta _{N,F} ^m$. 
\end{cor}

\begin{pf}
As in Corollary \ref{cor:CompMinSeq} we have that $\beta _{N,F} ^m< \infty$, then the conclusion follows from Theorem \ref{thm:CompLSCH}. 
\end{pf} 

\begin{rem}\label{rem:NonMinH}
As for the generalized second fundamental form, a priori, Corollary \ref{cor:CompMinSeq} does not ensure the existence of a minimizer since it can happen that the limit $m$-varifold $V$ is null. In  Section \ref{Sec:minH} we will see that, if $F(x,P,H^N)\geq C |H^N|^p$ for some $C>0$ and $p>m$, then this is not the case and we have a non trivial minimizer. 
\end{rem}

\section{Case $F(x,P,A) \geq C | A | ^p$ with $p>m$: non degeneracy of the minimizing sequence and existence of a $C^{1,\alpha}$ minimizer}\label{Sec:min}
Throughout this Section, $(\bar{N},g)$ stands for a compact $n$-dimensional Riemannian manifold isometrically embedded in some $\Rs$ (by Nash Embedding Theorem) and $N\subset \subset \bar{N}$ is a compact subset with non empty interior (as subset of $N$). Fix $m\leq n-1$; we will focus our attention and specialize the previous techniques to the case  
\begin{eqnarray}
F&:&G_m(N)\times \Rst\to \R^+  \text { is a function satisfying \eqref{def:F} } \nonumber\\
F(x,P,A)&\geq& C|A|^p \text{ for some } p>m \text{ and } C>0. \label{def:FAp}
\end{eqnarray}
Recall that we are considering the minimization problem
\begin{equation*}
\alpha_{N,F}^m:=\inf \left\{ \int_{G_m(N)} F(x,P,A) dV: V\in CV_m(N), V\neq 0 \text{ with generalized second fundamental form } A \right\}.
\end{equation*} 
Our goal is to prove the existence of a minimizer for $\alpha_{N,F}^m$, $F$ as in \eqref{def:FAp}.

Let  $\{V_k\}_{k\in \N} \subset CV_m(N)$ be a minimizing sequence of curvature varifolds  with generalized second fundamental forms $\{A_k\}_{k\in \N}$ such that 
$$\int_{G_m(N)} F(x,P,A_k) dV_k  \downarrow \a _{N,F} ^m;$$
from Corollary \ref{cor:CompMinSeq} we already know that there exists  $V \in CV_m(N)$ with generalized second fundamental form $A$ such that, up to subsequences,  

i) $(V_k, A_k) \rightharpoonup (V,A)$ in the weak sense of measure-function pairs,

ii) $\int_{G_m(N)} F(x,P,A) dV \leq \a _{N,F} ^m$. 

In order to have the existence of a minimizer we only have to check that $V$ is not the zero varifold; this will be done in the next Subsection \ref{SubSec:NonDeg} using the estimates of Section \ref{Sec:MF}.

\subsection{Non degeneracy properties of the minimizing sequence}\label{SubSec:NonDeg}
First of all, since $N\subset \Rs$, a curvature $m$-varifold $V$ of $N$ can be seen as a curvature varifold in $\Rs$ (for the precise value of the curvature function $B$ in $\Rs$ see Remark \ref{rem:BA}); as before we write $V=V(M,\th)$ where $M$ is a rectifiable set and $\th$ is the integer multiplicity function.  Let us call $H^{\R^s}$ the weak mean curvature of $V$ as integral $m$-varifold in $\Rs$ and, as in Section \ref{Sec:MF}, let us denote with $\mu=\mu_V={\cal H}^m\lfloor \theta=\pi_{\sharp} V$ the spatial measure associated to $V$ and with $\spt \mu$ its support. 

\begin{lem}\label{lem:Hbounded}
Let $N\subset \subset \bar{N}$ be a compact subset of the $n$-dimensional Riemannian manifold $(N,g)$ isometrically embedded in some $\Rs$ (by Nash Embedding Theorem) and fix $p>1$.  Then there exists a constant $C_{N,p}>0$ depending only on $p$ and $N$ such that for every $V=V(M,\th)\in CV_m(N)$  curvature $m$-varifold of $N$
$$\int_M \left|H^{\Rs}\right|^p d\mu\leq C_{N,p} \left( |V| + \int_{G_m(N)} |A|^p dV  \right).$$
\end{lem}

\begin{pf}
Recall (see Remark \ref{rem:BA}) that it is possible to write the curvature function $B$ of $V$  seen as curvature $m$-varifold of $\Rs$ in terms of the second fundamental form $A$ relative to $\bar{N}$ and the curvature of the manifold $\bar{N}$ seen as submanifold of $\Rs$ (the terms involving derivatives of $Q$):
\begin{equation*}
B_{ijk}=A^k_{ij}+A^j_{ik}+ P_{jl} P_{iq} \frac{\partial Q_{lk}}{\partial x_q}(x)+P_{kl} P_{iq} \frac{\partial Q_{lj}}{\partial x_q}(x).
\end{equation*}
From Remark \ref{rem:AH} the weak mean curvature $H^{\Rs}$, which is a vector of $\Rs$, can be written in terms of $B$ as

$$\left(H^{\Rs} \right)_i=\sum_{j=1}^S B_{jij}= \sum_{j=1}^S \left( A^j_{ji}+A^i_{jj}+ P_{il} P_{jq} \frac{\partial Q_{lj}}{\partial x_q}(x)+P_{jl} P_{jq} \frac{\partial Q_{li}}{\partial x_q}(x) \right) \quad i=1\ldots, S.$$
Notice that, since  $N\subset \subset \bar{N}$ is a compact subset of the manifold $\bar{N}$  smoothly embedded in $\Rs$, the functions $\frac{\partial Q_{lj}}{\partial x_m}$ are uniformly bounded by a constant $C_N$ depending on the embedding $N \hookrightarrow \Rs$; moreover the  $P_{jm}$ are projection matrices so they are also uniformly bounded and we can say that 

$$\left|\left( \sum_{j,l,m=1}^S P_{il} P_{jq} \frac{\partial Q_{lj}}{\partial x_q}+P_{jl} P_{jq} \frac{\partial Q_{li}}{\partial x_q}\right)_{i=1, \ldots, S} \right| \leq C_N$$
as vector of $\Rs$.
\\About the first term observe that, from the triangle inequality applied to the $\Rs$-vectors $(A^j_{ji})_{i=1,\ldots,S}$ ($j$ is fixed for each single vector),
$$\left|\left(\sum_{j=1}^S  A^j_{ji}\right)_{i=1,\ldots,S}\right|\leq \sum_{j=1}^S \left| (A^j_{ji})_{i=1,\ldots,S}\right|\leq S|A|$$
where, of course $|A|:= \sqrt{\sum_{i,j,k=1}^S (A^i_{jk})^2}\geq | (A^j_{ji})_{i=1,\ldots,S}|$ for all $j=1,\ldots, S$. The second adding term is analogous.
\\Putting together the two last estimates, by a triangle inequality, we have
$$\left| H^{\Rs} \right|\leq 2S |A|+ C_{N}.$$
Using the standard inequality $(a+b)^p\leq 2^{p-1} (a^p+b^p)$ for $a,b\geq0$ and $p>1$ given by the convexity of the function $t\mapsto t^p$ for $t\geq0,p>1$ we can write
\begin{equation}
\left| H^{\Rs} \right|^p \leq C_{N,p} \left( |A|^p+ 1\right).
\end{equation}
With an integration we get the conclusion.
\end{pf}

Using the estimates of Section \ref{Sec:MF} and the last Lemma we have uniform lower bounds on the mass and on the diameter of the spatial support of a curvature $m$-varifold $V\in CV_m(N)$ of $N$ with bounded $\int_{G_m(N)} |A|^p dV$, $p>m$.

\begin{thm}\label{thm:V>}
Let $N\subset \subset \bar{N}$ be a compact subset of the $n$-dimensional Riemannian manifold $(N,g)$ isometrically embedded in some $\Rs$ (by Nash Embedding Theorem) and fix $m\leq n-1$, $p>m$.

Then there exists a constant $C_{N,p,m}>0$ depending only on $p,m$ and on the embedding of $N$ into $\Rs$ such that $C_{N,p,m}\uparrow +\infty$ as $p\downarrow m$ and such that  for every $V=V(M,\th)\in CV_m(N)$  curvature $m$-varifold of $N$ with spatial measure $\mu$
\begin{eqnarray}
&i)&  \quad \diam_{\bar{N}} (\spt \mu) \geq \frac{1}{C_{N,p,m}\left( |V| + \int_{G_m(N)} |A|^p dV  \right)^{\frac{1}{p-m}}} \label{eq:diamM>}\\
&& \text{ where } \diam_{\bar{N}} (\spt \mu) \text{ is the diameter of } \spt \mu \text{ as a subset of the Riemannian manifold } \bar{N};\nonumber \\
&ii)& \quad C_{N,p,m} \; |V| \; \left( |V| + \int_{G_m(N)} |A|^p dV  \right)^{\frac{m}{p-m}}\geq 1 \label{eq:VA}. 
\end{eqnarray}
Notice that $\;ii)$ implies the existence of a constant $a_{N,m,p,\int|A|^p}>0$  depending only on $p,m$,  on $\int_{G_m(N)} |A|^p dV$ and on the embedding of  $N$ into $\Rs$, with  $a_{N,p,m,\int|A|^p} \downarrow 0$ if $p\downarrow m$ or if $\int_{G_m(N)} |A|^p dV\uparrow +\infty$   such that 
$$|V|\geq a_{N,p,m,\int|A|^p}>0.$$ 
\end{thm} 

\begin{pf}

$i)$ From Lemma \ref{lem:d>H}
\begin{equation*}
diam_{\bar{N}} (\spt \mu)\geq diam_{\Rs}(\spt \mu)\geq \frac{1}{C_{p,m}\Big(\int_M |H|^p d\mu \Big)^{\frac{1}{p-m}}}
\end{equation*}
where $C_{p,m}>0$ is a positive constant depending on $p,m$ and such that $C_{p,m}\to \infty$ if $p\downarrow m$. The conclusion follows plugging into the last inequality the estimate of Lemma \ref{lem:Hbounded}.

$ii)$ From Lemma \ref{A>H},
\begin{equation*}
|V| \geq  \frac{1}{C_{p,m}\Big(\int_M |H|^p d\mu \Big)^{\frac{m}{p-m}}}
\end{equation*}
with $C_{p,m}>0$ as above. The conclusion, again,  follows plugging into the last inequality the estimate of Lemma \ref{lem:Hbounded} and rearranging.
\end{pf}

\begin{cor}\label{cor:NonDeg}
Let $N\subset \subset \bar{N}$ be a compact subset with non empty interior, $int(N)\neq \emptyset$, of the $n$-dimensional Riemannian manifold $(N,g)$ isometrically embedded in some $\Rs$ (by Nash Embedding Theorem) and fix $m\leq n-1$.

Assume that the space $(N,g)$ does not contain any non zero $m$-varifold with null generalized second fundamental form and
consider a function $F:G_m(N)\times \Rst\to \R^+$   satisfying \eqref{def:F}, \eqref{def:FAp} and a corresponding minimizing sequence of curvature $m$-varifolds
$\{V_k\}_{k\in \N} \subset CV_m(N)$ with generalized second fundamental forms $\{A_k\}_{k\in \N}$ such that 
$$\int_{G_m(N)} F(x,P,A_k) dV_k  \downarrow \a _{N,F} ^m$$
( for the definition of $\a _{N,F} ^m$ see \eqref{eq:DefAlpha}). Then, called $\mu_k$ the spatial measures associated to $V_k$, there exists a constant $a_{N,F,m}>0$ depending only on $N$,$F$ and $m$ such that

\begin{eqnarray}
&i)&  \quad diam_{\bar{N}} (\spt \mu_k) \geq a_{N,F,m} \label{eq:diamMk>}\\
&ii)& \quad |V_k| \geq a_{N,F,m} \label{eq:Vk>}. 
\end{eqnarray}
\end{cor}

\begin{pf}
From Theorem \ref{thm:AreaEst} and the finiteness of $\a _{N,F} ^m$, since $(N,g)$ does not contain any non zero $m$-varifold with null generalized second fundamental form,
$$ |V_k| \leq C_{N,F,m} \int_{G_m(N)} F(x,P,A_k) dV_k \leq C_{N,F,m}$$
for some $C_{N,F,m}>0$ depending only on $N,F$ and $m$.
\\Moreover, since (by assumption \eqref{def:FAp}) $F(x,P,A)\geq C |A|^p$ for some $p>m$ and $C>0$, the boundness of $\int_{G_m(N)} F(x,P,A_k) dV_k$ implies that
$$ \int_{G_m(N)} |A_k|^p dV_k \leq C_{N,F,m} $$
for some $C_{N,F,m>}0$ depending only on $N,F$ and $m$.

The conclusion follows putting the last two inequalities into Theorem \ref{thm:V>}.
\end{pf}

\subsection{Existence and regularity of the minimizer}

Collecting Corollary \ref{cor:CompMinSeq} and Corollary \ref{cor:NonDeg} we can finally state and prove the first main Theorem \ref{thm:ExMin}.
\newline

\begin{pfnb} {\bf of Theorem \ref{thm:ExMin}}

If $a)$ is true we are done, so we can assume that $a)$ is not satisfied.

Let $\{V_k\}_{k\in \N} \subset CV_m(N)$  with generalized second fundamental forms $\{A_k\}_{k\in \N}$ be a minimizing sequence of $\a _{N,F} ^m$:  
$$\int_{G_m(N)} F(x,P,A_k) dV_k  \downarrow \a _{N,F} ^m.$$
Called $\mu_k$ the spatial measures associated to $V_k$, from Corollary  \ref{cor:NonDeg} we have the lower bounds:   
\begin{eqnarray}
&i)&  \quad diam_{\bar{N}} (\spt \mu_k) \geq a_{N,F,m} \nonumber\\
&ii)& \quad |V_k| \geq a_{N,F,m}, \nonumber 
\end{eqnarray}
for a constant $a_{N,F,m}>0$ depending only on $N$,$F$ and $m$.
Corollary \ref{cor:CompMinSeq} implies the existence of a curvature $m$-varifold $V=V(M,\th) \in CV_m(N)$ with generalized second fundamental form $A$ such that, up to subsequences,  

i) $(V_k, A_k) \rightharpoonup (V,A)$ in the weak sense of measure-function pairs of $N$,

ii) $\int_{G_m(N)} F(x,P,A) dV \leq \a _{N,F} ^m$.
\\The measure-function pair convergence implies the varifold convergence of $V_k \to V$ and the convergence of the associated spatial  measures 
$$\pi_\sharp V_k=:\mu_k \to \mu:=\pi_\sharp V \quad \text{weak convergence of Radon measures on $N$. } $$ 
It follows that
$$0<a_{N,F,m}\leq |V_k|=\mu_k(N)\to \mu(N)=|V|,$$
thus $V\neq 0$ is a minimizer for $\a_{N,F} ^m$. 

Proof of $b_1)$: $V$ is indecomposable in $CV_m(N)$. Otherwise $V=V_1+V_2$ with $0\neq V_1,V_2\in CV_m(N)$; the minimizing property of $V$ implies that, up to exchanging $V_1$ with $V_2$,  we must have  $\int_{G_m(N)} F(x,P,A) dV_1=\int_{G_m(N)} F(x,P,A) dV=  \a _{N,F} ^m$ and $V_2$ is a  non zero $m$-varifold with null generalized second fundamental form, contradicting the assumption that $a)$ does not hold. 

Proof of $b_2)$: since $N\hookrightarrow \Rs$ is properly embedded, the weak convergence $\mu_k \to \mu$ on $N$ implies the weak convergence of $\mu_k \to \mu$ as Radon measures on $\Rs$. From  mass bound on the $V_k$ and the bound on $\int_{G_m(N)}|A_k|^p dV_k$ given by the assumption \eqref{def:FAp} on $F$, Lemma \ref{lem:Hbounded} allows us to apply Proposition \ref{pro:HausConv} and we can say that the spatial supports
$$\spt \mu_k\to \spt \mu \quad \text{Hausdorff convergence as subsets of } \Rs.$$
Notice that, since $\bar{N}\hookrightarrow \Rs$ is embedded, the Hausdorff convergence of $M_k \to M$ as subsets of $\Rs$ implies 
$$\spt \mu_k\to \spt \mu \quad \text{Hausdorff convergence as subsets of } \bar{N},$$
and this implies that
$$0<a_{N,F,m}\leq \lim_k \diam_{\bar{N}} (\spt \mu_k) = diam_{\bar{N}} (\spt \mu).$$
hence $b2)$. 

Now the minimizer $V\in CV(N)$ is a non null curvature varifold on $N$ with generalized second fundamental form $A$ (relative to $\bar{N}$) in $L^p(V)$ for some $p>m$. Since $N\hookrightarrow \Rs$, $V$ can also be seen as a varifold in $\Rs$ and Remark \ref{rem:BA} tell that $V$ is actually a varifold with generalized curvature function $B$ given by
$$B_{ijk}=A^k_{ij}+A^j_{ik}+ P_{jl} P_{iq} \frac{\partial Q_{lk}}{\partial x_q}(x)+P_{kl} P_{iq} \frac{\partial Q_{lj}}{\partial x_q}(x)$$
where the terms of the type $ P_{jl} P_{iq} \frac{\partial Q_{lk}}{\partial x_q}(x)$ represent the extrinsic curvature
of $\bar{N}$ as a submanifold of $\Rs$ and, of course, are bounded on $N$ from the compactness:
$$\sup_{x\in N} \left| P_{jl} P_{iq} \frac{\partial Q_{lk}}{\partial x_q}(x)+P_{kl} P_{iq} \frac{\partial Q_{lj}}{\partial x_q}(x)\right| \leq C_N.$$
Hence, from triangle inequality,
$$|B|\leq 2 |A|+C_N$$
and  
$$|B|^p\leq C_{N,p} \;(|A|^p +1).$$
Using the mass bound $|V|=\lim_k |V_k|\leq C < \infty $, with an integration we get 
$$\int_{G_m(\Rs)} |B|^p dV < \infty.$$
Under this conditions Hutchinson shows in \cite{Hu2}  that $V$ is a locally a graph of multivalued $C^{1,\a}$ functions and that $b3)$ holds. 
\end{pfnb}

\section{Existence of an integral $m$-varifold with weak mean curvature minimizing $\int |H|^p$ for $p>m$}\label{Sec:minH}

As before, throughout this Section $(\bar{N},g)$ stands for a compact $n$-dimensional Riemannian manifold isometrically embedded in some $\Rs$ (by Nash Embedding Theorem) and $N\subset \subset \bar{N}$ is a compact subset with non empty interior (as subset of $N$). Fix $m\leq n-1$; analogously to Section \ref{Sec:min} we will focus our attention to the case  
\begin{eqnarray}
F&:&G_m(N)\times \Rs \to \R^+  \text { is a function satisfying \eqref{def:F} } \nonumber\\
F(x,P,H)&\geq& C|H|^p \text{ for some } p>m \text{ and } C>0. \label{def:FHp}
\end{eqnarray}
Recall that we are considering the minimization problem
\begin{equation*}
\beta_{N,F}^m:=\inf \left\{ \int_{G_m(N)} F(x,P,H^N) dV: V\in HV_m(N), V\neq 0 \text{ with weak mean curvature } H^N \text{ relative to } \bar{N} \right\}.
\end{equation*} 
Our goal is to prove the existence of a minimizer for $\beta_{N,F}^m$, $F$ as in \eqref{def:FHp}.

As in Section \ref{Sec:min} we consider a minimizing sequence $\{V_k\}_{k\in \N} \subset HV_m(N)$ of integral $m$-varifolds  with weak mean curvatures $\{H^N_k\}_{k\in \N}$ relative to $\bar{N}$ such that 
$$\int_{G_m(N)} F(x,P,H^N_k) dV_k  \downarrow \beta _{N,F} ^m;$$
from Corollary \ref{cor:CompMinSeqH} we already know that there exists  $V \in HV_m(N)$ with with weak mean curvature $H^N$ relative to $\bar{N}$ such that, up to subsequences,  

i) $(V_k, H^N_k) \rightharpoonup (V,H^N)$ in the weak sense of measure-function pairs,

ii) $\int_{G_m(N)} F(x,P,H^N) dV \leq \beta _{N,F} ^m$. 

In order to have the existence of a minimizer we only have to check that $V$ is not the zero varifold; this will be done analogously to Subsection \ref{SubSec:NonDeg} using the estimates of Section \ref{Sec:MF}.

As before, since $N\subset \Rs$, an integral $m$-varifold $V$ of $N$ with weak mean curvature $H^N$ relative to $\bar{N}$ can be seen as integral $m$-varifold  of $\Rs$ with weak mean curvature $H^{\Rs}$.  We write $V=V(M,\th)$ where $M$ is a rectifiable set and $\th$ is the integer multiplicity function; finally, as in Section \ref{Sec:MF}, let us denote with $\mu=\mu_V={\cal H}^m\lfloor \theta=\pi_{\sharp} V$ the spatial measure associated to $V$ and with $\spt \mu$ the spatial support of $V$. 

\begin{lem}\label{lem:HboundedH}
Let $N\subset \subset \bar{N}$ be a compact subset of the $n$-dimensional Riemannian manifold $(N,g)$ isometrically embedded in some $\Rs$ (by Nash Embedding Theorem) and fix $p>1$.  Then there exists a constant $C_{N,p}>0$ depending only on $p$ and $N$ such that for every $V=V(M,\th)\in HV_m(N)$  integral $m$-varifold of $N$ with weak mean curvature $H^N$ relative to $\bar{N}$ 
$$\int_M \left|H^{\Rs}\right|^p d\mu\leq C_{N,p} \left( |V| + \int_{G_m(N)} |H^N|^p dV  \right).$$
\end{lem}

\begin{pf}
By Definition \ref{def:VH} we can express
\begin{equation*}
(H^{\Rs})_i=(H^N)_i+ P_{jk} \frac{\partial Q_{ij}}{\partial x^k}
\end{equation*}
and from the triangle inequality
\begin{equation}\label{eq:TriIneH}
\left|H^{\Rs}\right|\leq \left|H^{N}\right|+ \left|P_{jk} \frac{\partial Q_{ij}}{\partial x^k}\right|;
\end{equation}
as vectors in $\Rs$. The second summand of the right hand side is a smooth function on the compact set $G_m(N)$ hence bounded by a constant $C_N$ depending on $N$:
$$\left|P_{jk} \frac{\partial Q_{ij}}{\partial x^k}\right|\leq C_N.$$
Using the standard inequality $(a+b)^p\leq 2^{p-1} (a^p+b^p)$ for $a,b\geq 0$ and $p>1$ we get
$$\left|H^{\Rs}\right|^p\leq C_{N,p} \big( 1+\left|H^{N}\right|^p \big) $$
which gives the thesis with an integration.
\end{pf}

\begin{rem}\label{rem:HV>}
An analogous result to Theorem \ref{thm:V>} holds, just replace $V=V(M,\th)\in CV_m(N)$ with $V=V(M,\th)\in HV_m(N)$ and $\int_{G_m(N)} |A|^p dV$ with $\int_{G_m(N)} |H^N|^p dV$.
\end{rem}

Now we can show the non degeneracy of the minimizing sequence for $\beta_{N,F}^m$, $F$ as in \ref{def:F}, \eqref{def:FHp}.

\begin{lem}\label{lem:NonDegH}
Let $N\subset \subset \bar{N}$ be a compact subset with non empty interior, $int(N)\neq \emptyset$, of the $n$-dimensional Riemannian manifold $(N,g)$ isometrically embedded in some $\Rs$ (by Nash Embedding Theorem) and fix $m\leq n-1$.

Assume that the space $(N,g)$ does not contain any non zero $m$-varifold with null weak mean curvature $H^N$ relative to $\bar{N}$ and
consider a function $F:G_m(N)\times \Rs\to \R^+$   satisfying \eqref{def:F}, \eqref{def:FHp} and a corresponding minimizing sequence of integral $m$-varifolds
$\{V_k\}_{k\in \N} \subset HV_m(N)$ with weak mean curvatures  $\{H^N_k\}_{k\in \N}$ relative to $\bar{N}$ such that 
$$\int_{G_m(N)} F(x,P,H^N_k) dV_k  \downarrow \beta _{N,F} ^m$$
( for the definition of $\beta _{N,F} ^m$ see \eqref{eq:DefBeta}). Then, called $\mu_k$ the spatial measures of $V_k$, there exists a constant $b_{N,F,m}>0$ depending only on $N$,$F$ and $m$ such that

\begin{eqnarray}
&i)&  \quad \diam_{\bar{N}} (\spt \mu_k) \geq b_{N,F,m} \label{eq:HdiamMk>}\\
&ii)& \quad |V_k| \geq b_{N,F,m} \label{eq:HVk>}. 
\end{eqnarray}
\end{lem}

\begin{pf}
From Theorem \ref{thm:AreaEstH} and the finiteness of $\beta _{N,F} ^m$, since $(N,g)$ does not contain any non zero $m$-varifold with null weak mean curvature $H^N$ relative to $\bar{N}$,
$$ |V_k| \leq C_{N,F,m} \int_{G_m(N)} F(x,P,H^N_k) dV_k \leq C_{N,F,m}$$
for some $C_{N,F,m}>0$ depending only on $N,F$ and $m$.
\\Moreover, since (by assumption \eqref{def:FHp}) $F(x,P,H^N)\geq C |H^N|^p$ for some $p>m$ and $C>0$, the boundness of $\int_{G_m(N)} F(x,P,H^N_k) dV_k$ implies that
$$ \int_{G_m(N)} |H^N_k|^p dV_k \leq C_{N,F,m} $$
for some $C_{N,F,m>}0$ depending only on $N,F$ and $m$.

The conclusion follows from the last two inequalities and Remark \ref{rem:HV>}.
\end{pf}

Now, collecting Corollary \ref{cor:CompMinSeqH} and Lemma \ref{lem:NonDegH}  we can finally state and prove Theorem \ref{thm:ExMinH}, namely the existence of a non trivial minimizer for $\beta_{N,F}^m$, $F$ as in \ref{def:F}, \eqref{def:FHp}.
\newline

\begin{pfnb} {\bf of Theorem \ref{thm:ExMinH}}

If $a)$ is true we are done, so we can assume that $a)$ is not satisfied.

Let $\{V_k\}_{k\in \N} \subset HV_m(N)$ with weak mean curvatures  $\{H^N_k\}_{k\in \N}$ be a minimizing sequence of $\beta _{N,F} ^m$:  
$$\int_{G_m(N)} F(x,P,H^N_k) dV_k  \downarrow \beta _{N,F} ^m.$$
Called $\mu_k$ the spatial measures of $V_k$, from Lemma  \ref{lem:NonDegH} we have the lower bounds:   
\begin{eqnarray}
&i)&  \quad diam_{\bar{N}} (\spt \mu_k) \geq b_{N,F,m}>0 \nonumber\\
&ii)& \quad |V_k| \geq b_{N,F,m}>0, \nonumber 
\end{eqnarray}
for a constant $b_{N,F,m}>0$ depending only on $N$,$F$ and $m$.
Corollary \ref{cor:CompMinSeqH} implies the existence of an integral $m$-varifold $V \in HV_m(N)$ with weak mean curvature $H^N$ relative to $\bar{N}$ such that, up to subsequences,  

i) $(V_k, H^N_k) \rightharpoonup (V,H^N)$ in the weak sense of measure-function pairs of $N$,

ii) $\int_{G_m(N)} F(x,P,H^N) dV \leq \beta _{N,F} ^m$.

Analogously to the proof of Theorem \ref {thm:ExMin}, one shows  that 
$$0<b_{N,F,m}\leq |V_k|=\mu_k(N)\to \mu(N)=|V|,$$
thus $V\neq 0$ is a minimizer for $\beta_{N,F} ^m$.
The proof of $b1)$ and $b2)$ are again analogous to the proof of the corresponding sentences in Theorem \ref{thm:ExMin}. Let us just comment on $b2)$: from the  mass bound on the $V_k$ and the bound on $\int_{G_m(N)}|H^N_k|^p dV_k$ given by the assumption \eqref{def:FHp} on $F$, Lemma \ref{lem:HboundedH} allows us to apply Proposition \ref{pro:HausConv} and, using the same tricks of Theorem \ref{thm:ExMin} we can say that the spatial supports
$$\spt  \mu_k\to \spt \mu \quad \text{Hausdorff convergence as subsets of } \bar{N},$$
and this implies that
$$0<b_{N,F,m}\leq \lim_k \diam_{\bar{N}} (\spt \mu_k) =\diam_{\bar{N}} (\spt \mu),$$
hence $b2)$. 
\end{pfnb}

\section{Examples and Remarks}\label{Sec:ExRem}
First of all let us point out that our setting includes, speaking about ambient manifolds, a large class of Riemannian manifolds with boundary. 
\begin{rem}\label{rem:ManBoundary}
Notice that if $N$ is a compact $n$-dimensional manifold with boundary then there exists an $n$-dimensional (a priori non compact) manifold  $\bar{N}$ without boundary such that $N$ is a compact subset of $\bar{N}$ (to define $\bar{N}$ just extend $N$ a little beyond $\partial N$ in the local boundary charts). Hence, given a compact $n$-dimensional Riemannian manifold $(N,g)$ with boundary  such that the metric $g$ can be extended in a smooth and non degenerate way (i.e. $g$ positive definite) up to the boundary $\partial N$, then $N$ is isometric to a compact subset of an $n$-dimensional Riemannian manifold $(\bar{N},\bar{g})$ without boundary. 

Thus all the Lemmas, Propositions and Theorems apply to the case in which the ambient space is a Riemannian manifold with boundary with the described non degeneracy property at $\partial N$.
\end{rem}

Now let us show that the main results Theorem \ref{thm:ExMin} and Theorem \ref{thm:ExMinH} are non empty, i.e we have examples of compact subsets of Riemannian manifolds where do not exist non zero varifolds with null weak mean curvature relative to $\bar{N}$ and a fortiori there exists no non zero varifold with null generalized second fundamental form. Let us start with an easy Lemma:

\begin{lem}\label{lem:NoH0NoA0}
Let $N\subset \subset \bar{N}$ be a compact subset of the $n$-dimensional Riemannian manifold $(N,g)$ isometrically embedded in some $\Rs$ (by Nash Embedding Theorem), fix $m\leq n-1$ and assume that $N$ contains no non zero $m$-varifold with null weak mean curvature relative to $\bar{N}$. Then $N$ does not contain any non zero $m$-varifold with null generalized second fundamental form.
\end{lem}   

\begin{pf}
We show that if the varifold $V$ has null generalized second fundamental form relative to $\bar{N}$ then $V$ also has null weak mean curvature relative to $\bar{N}$. Indeed let $V$ be a varifold on $N$ with generalized curvature function $B$ and second fundamental form $A$ relative to $\bar{N}$, then, from Remark \ref{rem:BA}, 
\begin{equation*}
B_{ijk}=A^k_{ij}+A^j_{ik}+ P_{jl} P_{ip} \frac{\partial Q_{lk}}{\partial x_p}(x)+P_{kl} P_{ip} \frac{\partial Q_{lj}}{\partial x_p}(x)
\end{equation*} 
where $P$ and $Q(x)$ are the projection matrices on $P\in G_m(N)$ and $T_x\bar{N}$.
Moreover, from Remark \ref{rem:AH}, $V$ has weak mean curvature as a varifold in $\Rs$ 
$$(H^{\Rs})_i=B_{jij};$$
hence, if the generalized second fundamental form $A$ is null, then
$$ (H^{\Rs})_i= P_{il} P_{jk} \frac{\partial Q_{lj}}{\partial x_k}(x)+P_{jl} P_{jk} \frac{\partial Q_{li}}{\partial x_k}(x). $$
It is not hard to check that the first summand of the right hand side is null (fix a point $x$ of $\bar{N}$ and choose a base of $T_x\bar{N}$ in which the Christoffel symbols of $\bar{N}$ vanish at $x$; write down the orthogonal projection matrix $Q$ with respect to this base and check the condition in this frame).
Thus $H^{\Rs}_i = P_{jk} \frac{\partial Q_{ij}}{\partial x^k}$ and Definition \ref{def:VH} gives
$$(H^N)_i=(H^{\Rs})_i - P_{jk} \frac{\partial Q_{ij}}{\partial x_k}(x) =0. $$
\end{pf}

Collecting Lemma \ref{lem:NoH0NoA0} and Remark \ref{rem:dVH} we can say that if a compact subset $N\subset \subset \bar{N}$  has a non zero $m$-varifold with null generalized second fundamental form, then a fortiori $N$ contains a non zero $m$-varifold with null weak mean curvature relative to $\bar{N}$, then a fortiori $N$ contains a non zero $m$-varifold with null first variation relative to $\bar{N}$ (recall, see Remark \ref{rem:StatVar}, that a varifold with null first variation is also called stationary varifold). Hence it is enough to give examples of compact subsets of Riemannian manifolds which do not contain any non zero $m$-varifold with null first variation relative to $\bar{N}$.

First, we mention two examples given by White in \cite{Whi} (for the proofs we refer to the original paper) next we will propose a couple of new examples which can be seen as a sort of generalization of White's ones. Recall that if $N$ is a compact Riemannian manifold with smooth boundary, $N$ is said to be \emph{mean convex} provided that the mean curvature vector at each point of $\partial N$ is an nonnegative multiple of the inward-pointing unit normal.

\begin{ex}\label{ex:Ric>}
Suppose that $N$ is a compact, connected, mean convex Riemannian manifold with smooth, nonempty  boundary, and that no component of $\partial N$ is a minimal surface. Suppose also that the dimension $n$ of $N$ is at most $7$ and that the Ricci curvature of $N$ is everywhere positive. Then $N$ contains no non zero $n-1$-varifold with null first variation relative to $N$ (i.e. stationary $n-1$-varifold).

More generally, if $N$ has nonnegative Ricci curvature, then the same conclusion holds unless $N$ contains a closed, embedded, totally geodesic hypersurface $M$ such that $Ric(\nu,\nu)=0$ for every unit normal $\nu$ to $M$ (where $Ric$ is the Ricci tensor of $N$).
\end{ex}

Minimal surfaces in ambient manifolds of the form $M\times \R$ have been deeply studied in recent years (see for example \cite{MeRo04}, \cite{MeRo05} and \cite{NeRo02}); notice that $M \times \R$ is foliated by the minimal surfaces $M\times \{z\}$.
In the second example we can see that very general compact subsets of ambient spaces admitting such foliations do not contain non zero  codimension 1 varifolds with null first variation.

\begin{ex}\label{ex:MxR}
Let $\bar{N}$ be an $n$-dimensional Riemannian manifold. Let $f:\bar{N} \to \R$ be a smooth function with nowhere vanishing gradient such that the level sets of $f$ are minimal hypersurfaces or, more generally, such that the sublevel sets $\{x: f(x) \leq z \}$ are mean convex. Let $N$ be a compact subset of $\bar{N}$ such that for each $z\in \R$, no connected component of $f^{-1}(z)$  is a minimal hypersurface lying entirely in $N$. Then $N$ contains no non zero  $n-1$-varifold with null first variation relative to $\bar{N}$.
 \end{ex}

Observe that both examples concern the non-existence of \emph{codimension 1} stationary varifolds: next we propose a couple of new examples in higher codimension. We need the following maximum principle for stationary (i.e. with null first variation) varifolds given by White, for the proof see \cite{Whi2}, Theorem 1. Before stating it recall that if $N$ is an $n$-dimensional Riemannian manifold with boundary $\partial N$, $N$ is said \emph{strongly $m$-convex} at a point $p\in \partial N$ provided 
$$k_1+k_2+\ldots+ k_m>0$$
where $k_1\leq k_2\leq\ldots \leq k_{n-1}$ are the principal curvatures of $\partial N$ at $p$ with respect to the unit normal $\nu_N$ that points into $N$.

\begin{thm}\label{thm:MaxPrinc}
Let $\bar{N}$ be a smooth Riemannian manifold of dimension $n$, let $N\subset \bar{N}$ be a smooth Riemannian $n$-dimensional manifold with boundary, and assume $p$ to be a point in $\partial N$ at which $N$ is strongly $m$-convex. Then $p$ is not contained in the support of any $m$-varifold in $N$ with null first variation relative to $\bar{N}$.
\end{thm}   
Actually the Theorem of White is more general and precise, but for our purposes this weaker version is sufficient.

Now are ready to state and prove the two examples.
\begin{thm}\label{thm:ExNRs}
Let $\bar{N}$ be an $n$-dimensional Riemannian manifold and consider as ambient manifold $\bar{N}\times \Rs, s>1$ with the product metric. Then any compact subset $N\subset \subset \bar{N}\times \Rs$ does not contain any non null stationary $n+k$-varifold,  $k=1,\ldots, s-1$  (i.e. $n+k$-varifold with null first variation relative to $\bar{N}\times \Rs$).
\end{thm}

\begin{pf}
Assume by contradiction that $V$ is a non null $n+k$-varifold in $N$ with null first variation in $\bar{N}\times \Rs$ for some $1\leq k\leq s-1$. Consider the function $\rh:\bar{N}\times \Rs \to \R^+$ defined as 
$$\bar{N}\times \Rs\ni (x,y)\mapsto \rh(x,y):=|y|_{\Rs}$$ 
where of course $|y|_{\Rs}$ is the norm of $y$ as vector of $\Rs$.
With abuse of notation, call $M\subset N$ the spatial support of $V$ (now $M$ may not be rectifiable, it is just compact); observe that, since $M$ is compact, then the function $\rh$ restricted to $M$ has a maximum $r>0$ at the point $(x_0,y_0)\in M\subset \bar{N}\times \Rs$ (observe that the maximum $r$ has to be strictly positive otherwise we would have a non null $n+k$-varifold in an $n$-dimensional space, which clearly is not possible by the very definition of varifold).
It follows that, called $\bar{N}_r$ the tube of center $\bar{N}$ and radius $r$
$$\bar{N}_r:=\{(x,y):x\in \bar{N}, |y|_{\Rs}\leq r\},$$
the spatial support of $V$ is contained in $\bar{N}_r$:
\begin{equation}\label{eq:MNr}
M\subset \bar{N}_r.
\end{equation}
Moreover $M$ is tangent to the hypersurface $C_r:=\partial \bar{N}_r=\{(x,y):x\in \bar{N}, |y|_{\Rs}= r\}$ at the point $(x_0,y_0)$. Observe that $C_r$ is diffeomorphic to $\bar{N}\times rS^{s-1}_{\Rs}$, where of course $rS^{s-1}_{\Rs}$ is the $s-1$-dimensional sphere of $\Rs$ of radius $r$ centered in the origin.

Using normal coordinates in $\bar{N}\times \Rs$ it is a simple exercise to observe that the principal curvatures of $C_r$ with respect to the inward pointing unit normal are constantly 
$$k_1=k_2=\ldots=k_n=0,\; k_{n+1}=k_{n+2}=\ldots=k_{s-1}=\frac{1}{r}$$
(just observe that the inward unit normal is $-\Theta$, where $\Theta$ is the radial vector which parametrizes $S^{s-1}_{\Rs}$; of course $-\Theta$ is constant respect to the $x$ coordinates; using normal coordinates one checks that the second fundamental form is made of two blocks: the one corresponding to $\bar{N}$ is null and the other one coincides with the second fundamental form of $S^{s-1}_{\Rs}$ as hypersurface in $\Rs$).

It follows that $C_r=\partial \bar{N}_r$ is strongly $n+k$-convex in all of its points, for all $1\leq k\leq n-1$; but $V$ is a non null $n+k$-varifold in $\bar{N}_r$ with null first variation relative to $\bar{N}$ and tangent to $C_r$ at the point $(x_0,y_0)\in C_{r}\cap M$. Fact which contradicts the maximum principle, Theorem \ref{thm:MaxPrinc}.    
 \end{pf}
 
As a corollary we have an example in all the codimensions in $\Rs$:

\begin{thm}\label{thm:ExRs}
Let $N\subset \subset \Rs$ be a compact subset of $\Rs$, $s>1$. Then, for all $1\leq m\leq s-1$, $N$ contains no non zero $m$-varifold with null first variation relative to $\Rs$. 
\end{thm}

\begin{pf}
Just take $\bar{N}:=\{x\}$ in the previous example, Theorem \ref{thm:ExNRs}, and observe that $\{x\}\times \Rs$ is isometric to $\Rs$.

Otherwise argue by contradiction as in the proof of Theorem \ref{thm:ExNRs} and observe that the support of the non zero $m$-varifold with null first variation is contained in a ball of $\Rs$ and tangent to its boundary, namely a sphere. Of course the sphere is strongly $m$-convex; it follows a contradiction with the maximum principle, Theorem \ref{thm:MaxPrinc}. 

\end{pf}
\begin{rem}\label{rem:ExStatVar}
Recall that if the ambient $n$-dimensional Riemannian manifold $N$ is compact without boundary, then Almgren proved in \cite {Alm} that for every $1\leq m<n$ there exists an integral $m$-varifold with null first variation relative to $N$. Moreover, in the same setting of compact $N$ and $\partial N=\emptyset$, Schoen and Simon \cite{ShSim81}, using the work of Pitts \cite{Pit81}, proved that $N$ must contain a closed, embedded hypersurface with singular set of dimension at most $n-7$. Hence, the isoperimetric inequality Theorem \ref{thm:AreaEstH} fails for such $N$ and the Theorem \ref{thm:ExMinH} is trivially true. However, as written above, there are many interesting examples of ambient manifolds with boundary where the Theorem is non trivial.
\end{rem} 

\begin{rem} \label{rem:A=0}
It is known that the ambient Riemannian $n$-manifolds, $n \geq 3$ (with or without boundary) which contain a smooth $m$-dimensional submanifold, $m\geq2$, with null second fundamental form (i.e a totally geodesic submanifold) are quite rare. It could be interesting to show the same in the context of varifolds, that is to prove that the ambient compact Riemannian $n$-manifolds, $n \geq 3$ (with or without boundary) which contain a non zero (a priori non rectifiable) $m$-varifold, $m\geq2$, with null second fundamental form relative to $N$ (see Definition \ref{def:VC}) are quite rare. This fact would imply the existence of a larger class of ambient Riemannian manifolds where the isoperimetric inequality Theorem \ref{thm:AreaEst} holds and the main Theorem \ref{thm:ExMin} is non trivial.
\end{rem}

\end{document}